\documentclass[11pt,a4paper]{article}
\usepackage{amsmath, amsfonts, amssymb,amsthm}
\usepackage{a4wide}
\usepackage{parskip}
\usepackage{enumitem}
\usepackage{xcolor}
\usepackage{pstricks}
\usepackage{hyperref}
\usepackage{latexsym}
\usepackage{cite}

\newcommand{\e}{\epsilon}

\newcommand{\al} {\alpha}
\newcommand{\ba} {\beta}
\newcommand{\de} {\delta}
\newcommand{\ga} {\gamma}
\newcommand{\Ga} {\Gamma}

\newcommand{\De} {\Delta}
\newcommand{\la} {\lambda}

\newcommand{\noi} {\noindent}
\newcommand{\na} {\nabla}
\newcommand{\ds} {\displaystyle}

\newcommand{\RR}{{\mathbb R}}
\newcommand{\te} {\theta}
\newcommand{\NN}{{\mathbb N}}

\newtheorem{defi}{Definition}[section]
\newtheorem{thm}{Theorem}[section]
\newtheorem{rem}{Remark}[section]
\newtheorem{lem}{Lemma}[section]
\newtheorem{prop}{Proposition}[section]

\numberwithin{equation}{section}
\def\proof{\noindent{\textbf{Proof. }}}

\begin{document}
	\setlength{\abovedisplayskip}{3pt}
	\setlength{\belowdisplayskip}{3pt}
	\date{}
	\vspace{0.01in}
		\title{High energy solutions for $p$-Kirchhoff elliptic problems with Hardy-Littlewood-Sobolev nonlinearity}

		\author{ {\bf Divya Goel$\,^{1,}$\footnote{e-mail: {\tt divya.mat@iitbhu.ac.in}},  Sushmita Rawat $\,^{2,}$\footnote{e-mail: {\tt sushmita.rawat1994@gmail.com}} and 
				K. Sreenadh$\,^{2,}$\footnote{e-mail: {\tt sreenadh@maths.iitd.ac.in}}} \\ $^1\,$Department of Mathematical Sciences, Indian Institute of Technology (BHU),\\ Varanasi 221005, India\\
			$^2\,$Department of Mathematics, Indian Institute of Technology Delhi,\\ Hauz Khas, New Delhi 110016, India.}

		\maketitle
		\begin{abstract}
			This article deals with the study of the following Kirchhoff-Choquard problem:
			\begin{equation*}
				\begin{array}{cc}
					\displaystyle  M\left(\, \int\limits_{\mathbb{R}^N}|\nabla u|^p\right) (-\Delta_p) u + V(x)|u|^{p-2}u = \left(\, \int\limits_{\RR^N}\frac{F(u)(y)}{|x-y|^{\mu}}\,dy \right) f(u), \;\;\text{in} \; \mathbb{R}^N,\\
					u > 0, \;\; \text{in} \; \mathbb{R}^N,
				\end{array}
			\end{equation*}
			where $M$ models Kirchhoff-type nonlinear term of the form $M(t) = a + bt^{\theta-1}$, where $a, b > 0$ are given constants; $1<p<N$, $\Delta_p = \text{div}(|\nabla u|^{p-2}\nabla u)$ is the $p$-Laplacian operator; potential $V \in C^2(\mathbb{R}^N)$; $f$ is monotonic function with suitable growth conditions. We obtain the existence of a positive high energy solution for $\theta \in \left[1, \frac{2N-\mu}{N-p}\right) $ via the Poho\v{z}aev manifold and linking theorem.
  Apart from this, we also studied the radial symmetry of solutions of the associated limit problem. 
			
			\medskip

			
			\medskip
			
			
		\end{abstract}
		\newpage
		\section{Introduction}

		The classical work by W. Ding and W. M. Ni in \cite{ding} and P. Rabinowitz in \cite{rabino}   opens a new stage  in solving partial differential equations. The authors proved the existence of critical points  using constrained minimization. The innovative idea was to establish that the mountain pass min-max level of the functional $J$, defined on the Hilbert space $H$ and 
		associated with the equation, is equal to the minimum of $J$
		restricted to the so-called Nehari manifold $\mathcal{N}= \{u \in H|  \langle J^\prime (u), u\rangle =0 \}$. Further, this technique has solved numerous problems, given the homogeneous and superlinear condition at infinity of the  nonlinearity. 
		
		In contrast, if the  nonlinear term is not homogeneous and superlinear,  then we cannot employ the strategy of minimizing over the Nehari manifold. To overcome this, Jeanjean and Tanaka \cite{tanaka} defined  the Poho\v{z}aev manifold and obtained solutions as minimizers of the functional over this manifold. Afterward, many authors used this technique to prove the existence of ground-state solutions to elliptic equations. In this article, we are interested in studying the equations  for  the $p$-Kirchhoff type operator with  nonlocal  non-homogeneous nonlinearity  and without the superlinear condition at infinity. 
 
 Elliptic equations of the Kirchhoff type have been studied extensively in the literature. Consider the following model problem
		\begin{equation}\label{Eeq1.3}
			-\left(  a +b\int\limits_{\RR^N}|\nabla u|^2\,dx  \right)\Delta u + V(x) u = f(x, u) \text{ in } \RR^N
		\end{equation}
		where $f(x, t)$ is subcritical and superlinear at $t=0$. For this class of problems, extensive literature exists on the existence and multiplicity of results using variational methods. For instance, one can see \cite{wu, Chen_li, sun_tang, tang_cheng, naimen, he_zou, alves, rawat_sreenadh, goel_rawat_sreenadh} and the references therein. 
		In \eqref{Eeq1.3},  Kirchhoff term is homogeneous with degree 4, so   to obtain the geometric structure and the boundedness of Palais-Smale sequences for the energy functional, one usually assumes one of the following conditions
		\begin{enumerate}
			\item[(SL)] \text{4-superlinear condition:}
			$$	\lim\limits_{|t| \to \infty} \frac{F(x,t)}{t^4} = \infty \;\text{uniformly in}\; x\in \RR^N,\; \text{where}\; F(x,t)=\int_{0}^{t}f(x,s)ds;$$
			\item[(AR)] \text{Ambrosetti–Rabinowitz type condition:} For all $t \in \RR$ there exists $\la >4$ such that 
			\begin{equation*}
				f(x,t)t \geq \ \la F(x,t) >0;
			\end{equation*}
			\item[(VC)] \text{Variant convex condition:}
			$\frac{f(x,t)}{|t|^3} \; \text{is strictly increasing for}\; t\in \RR\backslash\{0\}.$
		\end{enumerate}
		Even though these conditions seem to be fundamental, they are restrictive in nature. For example, $f(x,t)=|t|^{q-2}t$, when $2< q \leq 4$, the above conditions do not hold. So the elementary question arises as to whether one can obtain a ground-state solution without these conditions.

In 2014, Li and Ye \cite{li_ye_subcritical} worked in this direction for the Laplacian and  obtained a positive ground state solution of \eqref{Eeq1.3} on $H^1(\RR^3)$ for  $f(x,u) = |u|^{q-2}u$ and $3<q \leq 4$ with  some suitable assumptions on $V$. Authors used  the minimizing argument on a new manifold called Nehari-Poho\v{z}aev manifold, defined as the set of all functions $u\in H^1(\RR^3)$ satisfying $ \langle I'(u),u \rangle + Q(u)= 0,$
	where $I(u)$ is the energy functional and $Q(u):= \lim\limits_{n \to \infty}\langle I'(u), \phi_n(x\cdot \nabla u)  \rangle$ is the Poho\v{z}aev functional for \eqref{Eeq1.3} and $\{\phi_n\}$ are suitable test functions. 
  Subsequently,  Guo \cite{guo} generalized the results of \cite{li_ye_subcritical}  for a general nonlinearity $f(u)$ with  the following assumptions: 
 \begin{equation*}
		f\in C^1 \text{ and }	 \frac{f(s)}{s} \;\text{is strictly increasing in} \; (0,\infty). 
	\end{equation*}
 Here, the author  established that \eqref{Eeq1.3} admits a least energy solution using Jeanjean's \cite{jeanjean-mon} monotonicity approach on a Nehari-Poho\v{z}aev type manifold.
 %
\noi However, in \cite{hongyu}, Ye, with the superlinear condition at infinity, observed  that the groundstate solution for \eqref{Eeq1.3} does not exist. Hence the author obtained a positive high energy solution by assuming suitable conditions for the potential function $V$. This was achieved by using the linking theorem on Nehari-Poho\v{z}aev type manifold.
For more results, one can refer to \cite{li, azzollini,tang_chen, he_zou}.

	On the other hand, problems that involve two nonlocal terms have been studied by L\"u \cite{lu_hartree_NL} with the potential  $ V(x) = 1+\la g(x)$. More precisely, he considered the problem 
	\begin{equation*}
		-M\left(\,\int\limits_{\RR^3} |\nabla u|^2\,dx\right)\De u +V(x)u= \left(\,\int\limits_{\RR^3}\frac{|u(y)|^q}{|x-y|^{3 -\al}}\,dy\right)|u|^{q-2}u, \; u>0\; \text{in}\;\RR^3
	\end{equation*}
	where $M(t) = a +bt$, $q\in (2,3+ \al)$, $\la >0$ is a parameter and $g(x)$ is a non-negative steep potential well function. 
	For  $\la$ large enough, L\"u proved the existence of ground state solutions using the Nehari manifold and the concentration-compactness principle. Since the associated functional is not bounded below, the assertions used here can not be carried forward for the case of $q \in (1	+ \al/3,	2]$. But using Nehari-Poho\v{z}aev manifold with some strong  assumptions on the potential $V(x)$, Chen \&  Liu\cite{chen_liu},  obtained a ground state solution in $\RR^3$ for the complete range, $q \in (1 + \al/3, 3 + \al)$. For more recent results on ground state solutions, refer to \cite{chen_zhang_tang,zhang_xie_jiang}. 

	Inspired  by the works described above, 
	in this paper, we study the existence of high energy solutions for a class of $p$-Laplace equations with a general Kirchoff terms $M(t) = a+ bt^{\te-1}$ and the nonlinearity that does not satisfy (SL), (AR), and (VC) assumptions.
	%
 Precisely, we study  the following class of Kirchhoff problems with the general Choquard term:
\begin{equation}\label{E1.1}    
		\begin{array}{rllll}\ds -M\left( \int_{\mathbb R^N} |\na u|^p dx \right)  \Delta_p u + V(x)|u|^{p-2}u &=  \left(\,\displaystyle\int\limits_{\mathbb{R}^N} \frac{F(u)(y)}{|x-y|^\mu}\,dy \right) f(u),\,  \; u>0 \;\text{in} \; \mathbb{R}^{N},
		\end{array}
	\end{equation}
	where $N > \max\{2, p\}$ with $1<p <N$, $a, b, \theta$ are positive parameters and $0 < \mu < N$. Here,  $\Delta_p$ is the $p$-Laplacian operator, defined as $\Delta_pu =$ div$(|\nabla u|^{p-2}\nabla u)$. The function $f\in C(\RR , \RR) $ is  such that $f(t) \equiv 0$ for $t \leq 0$ and $F$ is the primitive of $f$. The general function $f$ satisfies the following growth conditions:
	\begin{enumerate}
		\item [$\mathbf{(\mathcal{F}_1)}$] $\ds\lim\limits_{t \to 0}\frac{f(t)}{t^{p_{*, \mu}-1}} = 0$;
		\item [$\mathbf{(\mathcal{F}_2)}$] $f$ is a monotonic function and for some positive constant $C$ and $q\in \left( p_{*, \mu}, p^*_\mu\right)$ $$\ds|f(t)t| \leq C |t|^{p_{*, \mu}} + C |t|^q;$$
	\end{enumerate}
	where $ p_{*, \mu}= \frac{p(2N-\mu)}{2N}$ is the lower critical exponent and $p^{*}_{\mu} = \frac{p(2N-\mu)}{2(N-p)}$ is the upper critical exponent in the sense of the Hardy–Littlewood–Sobolev inequality.
	Moreover, we assume that the potential function $ V \in C^2(\RR^N, \RR)$ verifies:
	\begin{enumerate}
		\item [$\mathbf{(\mathcal{V}_1)}$]  $\lim\limits_{|x| \to \infty}V(x) := V_\infty > 0$;
		\item [$\mathbf{(\mathcal{V}_2)}$]  $\nabla V(x)\cdot x \leq 0$, for all $x\in \RR^N$ with the strict inequality on a subset of positive Lebesgue measure;
		\item [$\mathbf{(\mathcal{V}_3)}$] $NV(x) + \nabla V(x)\cdot x \geq NV_\infty$, for all $x \in \RR^N$;
		\item [$\mathbf{(\mathcal{V}_4)}$] $ N\nabla V(x)\cdot x + x\cdot H(x)\cdot x \leq 0$, for all $x \in \RR^N$ where $H$ represents the Hessian matrix of the function $V$.
	\end{enumerate}	
	This paper presents a new approach to show  the existence results to \eqref{E1.1} for all $N\geq 3$. 
	We established the existence of high energy solutions to \eqref{E1.1}
 for a larger class of nonlinearity $f$. 
	Most of the literature mentioned above for Kirchhoff problems is for the case when $N =3$, as the case $N\geq 3$ is more intricate. The techniques and ideas presented in \cite{hongyu} cannot be extended to higher dimensions using the Nehari-Poho\v{z}aev type manifold. This is primarily because the super quadratic growth in the Kirchhoff operator dominates the nonlinearity. 
 %
	 We study the problem with the Kirchhoff term
	$M(t) = a+ bt^{\te-1},$
	for all values of  $\te \in \left[1 , \frac{2N-\mu}{N-p}\right)$.
	We remark that we incorporated the case $\te =2$, that is, $M(t) = a+ bt$ when $\mu < \min\{N, 2p\}$. 
	To successfully establish the existence of a solution to \eqref{E1.1}, we need first to study the following  limit problem 
	\begin{equation}  \label{E1.2}
		\ds - M\left(\int\limits_{\RR^N}|\nabla u|^p\right) \Delta_p u + V_\infty|u|^{p-2}u =  \left(\,\,\int\limits_{\mathbb{R}^N} \frac{F(u)(y)}{|x-y|^\mu}\,dy \right) f(u)\,  \;\; \text{in} \; \mathbb{R}^{N}.
	\end{equation}
	 We study the existence of a solution to  this limit   problem using the concentration compactness arguments of Lions and the radial symmetry nature of the solution. 
  Proving the radial symmetry nature of the solutions for this class of equations is an independent topic of interest.   In the Laplacian case ( $p=2, M=1$),  the radial nature of the solution was established long back by Li \& Ni \cite{liandni}  for general ground state solutions, under the condition $f^\prime(s) <0$. 
	However, the same approach cannot be carried forward for the $p$-Laplacian case ($p\not =2, M=1$).  In \cite{badiale}, the radial symmetry of solutions for a class of $p$-Laplacian equations was proved under the assumption that  solution $v$ has only one critical point in  $ \RR^N$. Subsequently, many researchers obtained radial symmetry result by using the $a$priori estimates on the solutions. But  in our case, we establish the radial nature of the solution  by aptly using the  non-local nature of  nonlinearity. 
	 Our arguments are based on
	the Euler-Lagrange equation satisfied by polarizations and   equality cases in polarization inequalities. We also establish some inequalities in the frame of polarization which we didn't find explicitly in the literature.  We dedicate Section 3 and Section 4 to establish that the ground state solution of \eqref{E1.2} is radially decreasing. We develop a unified approach in this paper to study both \eqref{E1.1} and \eqref{E1.2}. In the next section, we give some preliminary framework of the problem  and state the main results of the article. 

	\section{Preliminaries and Main result}
	This section of the article is intended to provide the variational  setting. Refer to \cite{leib, zhu_cao} for the complete and rigid details. Further in this section, we state the main results
	of the current article with a short sketch of the proof.  The functional space associated to this problem is $$W^{1,p}(\RR^N) = \{u \in L^P(\RR^N): \int\limits_{\RR^N}|\nabla u|^p < \infty \}$$ with the norm,
	\begin{equation*}
		\|u\|^p = \ds a\int\limits_{\RR^{N}} |\nabla u|^p\,dx + \int\limits_{\RR^{N}} V(x)| u|^p\,dx.
	\end{equation*}
	It is easy to see from the condition $(\mathcal{V}_2)-(\mathcal{V}_3)$, $V_\infty \leq V(x)$ and from $(\mathcal{V}_1)$ we get $\sup\limits_{\RR^{N}}V(x) <\infty$, which infers that the norm is well-defined.
	\begin{prop}
		\textbf{(Hardy-Littlewood-Sobolev inequality)} Let $t$, $r > 1$ and $0 < \mu < N$ with $\frac{1}{t} + \frac{\mu}{N} + \frac{1}{r} = 2$, $f \in L^t(\mathbb{R}^N)$ and $h \in L^r(\mathbb{R}^N)$. Then there exists a sharp constant $C(t, r, \mu, N)$ independent of $f$, $h$ such that
		\begin{equation*}
			\iint\limits_{\mathbb{R}^{2N}} \dfrac{f(x)h(y)}{|x-y|^ \mu}\,dxdy \leq C(t, r, \mu, N)\|f\|_{L^t(\mathbb{R}^N)}\|h\|_{L^r(\mathbb{R}^N)}.
		\end{equation*}
	\end{prop}
	By Hardy-Littlewood-Sobolev inequality and condition $(\mathcal{F}_2)$, there exists constant $C>0$ such that for each $u \in W^{1,p}(\RR^N)$,
	\begin{equation}\label{Eeq2.1}
		\begin{aligned}
			\iint\limits_{\RR^{2N}} \frac{F(u)F(u)}{|x-y|^ \mu}\,dxdy 
			\leq & C \left[\;\int\limits_{\RR^N}|u|^p + |u|^{\frac{2Nq}{2N-\mu}} \right]^\frac{2\cdot p^*_{\mu}}{p^*}
			\leq &C\|u\|_{L^p}^{2\cdot p_{*,\mu}}+ C \|u\|_{L^{\frac{2Nq}{2N-\mu}}}^{2q}.
		\end{aligned}
	\end{equation}

	\begin{lem}\label{Elem3.2} Let $1 < p \leq \infty$, $1 \leq q < \infty$, with $q \neq \dfrac{Np}{N-p}$ if $p<N$. Assume that $v_n$ is bounded in $L^q (\RR^N )$, $|\nabla v_n|$ is bounded in $L^p(\RR^N )$ and
		$\ds\sup \limits_{y\in \RR^N}\int\limits_{B_R(y)}|v_n|^q\;dx \to 0$, for some $R > 0$ as $n \to \infty$.
		Then $v_n \to 0$ in $L^\al(\RR^N )$, for $\al \in (q, \frac{Np}{N-p}). $
	\end{lem}

	The energy functionals associated with the problems \eqref{E1.1}  and \eqref{E1.2} are  $J, J_\infty : W^{1,p}(\RR^N) \rightarrow \mathbb{R}$ defined as
	\begin{equation*}
		\begin{aligned}
			&	J(u) =  \frac{a}{p}\int\limits_{\RR^N} |\nabla u|^p dx + \frac{b}{p\te}\left(\, \int\limits_{\RR^N} |\nabla u|^p dx\right)^{\te} +\frac{1}{p}\int\limits_{\RR^N}V(x)|u|^pdx -\frac{1}{2}\iint\limits_{\RR^{2N}} \frac{F(u)F(u)}{|x-y|^ \mu}dxdy,\\
			& \ds J_{\infty}(u) = \frac{a}{p}\int\limits_{\RR^N} |\nabla u|^p dx + \frac{b}{p\te}\left(\, \int\limits_{\RR^N} |\nabla u|^p dx\right)^{\te} +\frac{V_\infty}{p}\int\limits_{\RR^N}|u|^pdx -\frac{1}{2}\iint\limits_{\RR^{2N}} \frac{F(u)F(u)}{|x-y|^ \mu}dxdy.
		\end{aligned}
	\end{equation*}
	Observe that both the functionals  $J, J_\infty \in C^1$, that is, for any $\phi \in C^{\infty}_0(\RR^N)$
		\begin{equation*}
			\begin{aligned}
				\ds \langle J'(u), \phi\rangle= & a\int\limits_{\RR^{N}}|\nabla u|^{p-2}\nabla u\cdot\nabla\phi\,dx + b\left(\, \int\limits_{\RR^N} |\nabla u|^p dx\right)^{\te-1}\int\limits_{\RR^{N}}|\nabla u|^{p-2}\nabla u\cdot\nabla\phi\,dx\\ 
				& + \int\limits_{\RR^{N}}V(x)|u|^{p-2}u\phi\,dx - \iint\limits_{\RR^{2N}}\frac{F(u(y))f(u(x))\phi(x)}{|x-y|^ \mu}dxdy,
			\end{aligned}
		\end{equation*}
		and 
		\begin{equation*}
			\begin{aligned}
				\ds \langle J_\infty'(u), \phi\rangle= & a\int\limits_{\RR^{N}}|\nabla u|^{p-2}\nabla u\cdot\nabla\phi\,dx + b\left(\, \int\limits_{\RR^N} |\nabla u|^p dx\right)^{\te-1}\int\limits_{\RR^{N}}|\nabla u|^{p-2}\nabla u\cdot\nabla\phi\,dx\\ 
				& + V_\infty\int\limits_{\RR^{N}}|u|^{p-2}u\phi\,dx - \iint\limits_{\RR^{2N}}\frac{F(u(y))f(u(x))\phi(x)}{|x-y|^ \mu}dxdy.
			\end{aligned}
	\end{equation*}
We also define the following Pohozaev manifolds associated with problem \eqref{E1.1}  and \eqref{E1.2}
$$\mathcal{P} = \{u \in W^{1,p}(\RR^N)\backslash\{0\}: P(u)= 0\} \;\; \text{and}\; \mathcal{P_\infty} = \{u \in W^{1,p}(\RR^N)\backslash\{0\}: P_\infty(u)= 0\}.$$
Here, $P(u)$ and $P_\infty(u)$ are the Poho\v{z}aev identities associated with problem \eqref{E1.1}  and \eqref{E1.2} respectively, (see \cite[Theorem 3]{Moroz1})
\begin{equation}
\begin{aligned}\label{Eeq1.5}
		P(u): = & \frac{(N-p)a}{p}\int\limits_{\RR^N} |\nabla u|^p dx + \frac{(N-p)b}{p}\left(\, \int\limits_{\RR^N} |\nabla u|^p dx\right)^{\te} +\frac{1}{p}\int\limits_{\RR^N}\left( N V(x)+ \nabla V(x)\cdot x\right) |u|^p\\
		& -\frac{2N-\mu}{2}\iint\limits_{\RR^{2N}} \frac{F(u)F(u)}{|x-y|^ \mu}dxdy,
	\end{aligned}
\end{equation} 
and
\begin{equation}
	\begin{aligned}\label{Eeq1.6}
		P_\infty(u): =&  \frac{(N-p)a}{p}\int\limits_{\RR^N} |\nabla u|^p dx + \frac{(N-p)b}{p}\left(\, \int\limits_{\RR^N} |\nabla u|^p dx\right)^{\te} +\frac{NV_\infty}{p}\int\limits_{\RR^N}|u|^pdx\\
		&-\frac{2N-\mu}{2}\iint\limits_{\RR^{2N}} \frac{F(u)F(u)}{|x-y|^ \mu}dxdy.
	\end{aligned}		
\end{equation}


With these preliminaries, we state our main results
\begin{thm}\label{Ethm1.1}
	Assume that $f$ and $V$ satisfy $(\mathcal{F}_1)-(\mathcal{F}_2)$ and $(\mathcal{V}_1)-(\mathcal{V}_4)$. Then, problem \eqref{E1.2} has a positive ground state solution that is radially symmetric.
\end{thm}
\begin{thm}\label{Ethm1.2}
	Assume that $f$ and $V$ satisfy $(\mathcal{F}_1)-(\mathcal{F}_2)$ and $(\mathcal{V}_1)-(\mathcal{V}_4)$ along with
	\begin{enumerate}
		\item [$\mathbf{(\mathcal{V}_5)}$] there exists a constant $T > 1$ such that
		\begin{equation*}
			\sup\limits_{x \in\RR^N} V(x) \leq V_\infty + \frac{p J_\infty(\overline{u})}{T^N \int\limits_{\RR^N}|\overline{u}|^p},
		\end{equation*}  
	\end{enumerate}
	where $\overline{u}$ is the ground state solution of the functional $J_\infty$ obtained in Theorem \ref{Ethm1.1}.
	Then there exists at least one positive solution to problem \eqref{E1.1}, which is a high energy solution.
\end{thm}
A standard methodology to obtain solutions using variational techniques    consists of looking for minimizers of the functional.  Here, using the  (SL) and (AR)  conditions, one can prove the boundedness of the minimizing sequences, which will subsequently give us solutions. But in the above-mentioned results, due to  the absence of conditions (SL) or (AR), verifying the mountain-pass geometry and obtaining the boundedness of the Palais-Smale sequence for the functionals $J$ and $J_\infty$ is a tedious job. Moreover, the lack of conditions (VC) and $(f_3)-(f_4)$  makes the  classical  method of  the Nehari manifold inappropriate. Taking into account all these points, we observed that $J_\infty$ satisfies the  mountain-pass geometry of the functional  with the scaling function $u_{t}$ defined as 
$$u_{t}(x):= u(t^{-1}x),$$
for $u \in W^{1,p}(\RR^N)\backslash\{0\}$ and $t>0$. 
Keeping this in mind, we take minimum over the Poho\v{z}aev manifold   $ \mathcal{P_\infty}.$
Using the assumptions on $f$ and $V_\infty$, we conclude that 
\begin{enumerate}
	\item[(1)]  $\mathcal{P_\infty}$ is a natural constraint, in the sense that $J_\infty$ is coercive and bounded below on $\mathcal{P_\infty}$
	\item[(2)] Every critical point of $J_\infty$ restricted to $\mathcal{P_\infty}$ is a critical point of $J_\infty$
	\item[(3)] $ \mathfrak{p_\infty}:= \inf\limits_{\mathcal{P_\infty}}J_\infty(u) = c_\infty$ where $c_\infty$ is the mountain pass level.  
\end{enumerate}
Further, to construct  a bounded Palais-Smale sequence for $J_\infty$, we employ Jeanjean's \cite{jeanjean} technique and prove the existence of a non-trivial critical point of functional $J_\infty$ by the concentration-compactness lemma. This motivates us to use the same type of approach for the existence of a solution to \eqref{E1.1}. That is, we look for the solution in the Poho\v{z}aev manifold $\mathcal{P}$ and
with the same type of scaling, one can establish that   the functional $J$ satisfies the mountain-pass geometry.
However, we encounter that the functional $J$  couldn't attain the mountain pass level, thus indicating the absence of a ground state solution for equation \eqref{E1.1}. To overcome this obstacle, we establish that $J$ satisfies the Palais-Smale condition above the min-max level, leading to  solutions with higher energy.  To address the issue of compactness, we employ the splitting lemma (refer to Lemma \ref{Elem5.9}). Additionally, we utilize the linking theorem in conjunction with a Barycenter mapping to achieve the goal in Theorem \ref{Ethm1.2}. To successfully apply the machinery of Barycenter mapping, we need the radial nature of ground state solution to \eqref{E1.2}. 
The salient feature of the approach is the existence of the radial solution to \eqref{E1.2}. To the best of our knowledge, there has been no attempt till now on the existence of high energy solutions to a generalized class of Kirchhoff-Choquard equation driven by the $p$-Laplacian operator.





\begin{rem} Conditions $(\mathcal{V}_1)$,$(\mathcal{V}_2)$ and $(\mathcal{V}_3)$ imply that
	$$\nabla V(x)\cdot x \to 0, \;\text{if} \;|x|\to \infty.$$
\end{rem}
\begin{rem}
    There are a number of functions that satisfy $(\mathcal{V}_1)-(\mathcal{V}_5)$. For example, $V(x)= V_\infty + \frac{C}{|x|^2+1}$, where $C = \frac{p J_\infty(\overline{u})}{T^N \int\limits_{\RR^N}|\overline{u}|^p}$.
\end{rem}

Throughout the paper, we make use of the following notations:
\begin{itemize}
	\item $u^{\pm} := \max\{\pm u, 0\}$ and $(H^*, \|\cdot\|_*)$ denotes the dual space of $(H, \|\cdot\|)$.
	\item For any $u \in W^{1,p}(\RR^N)$
	\begin{equation*}
		\|F(u)\|_{0}:= \iint\limits_{\RR^{2N}} \frac{F(u)F(u)}{|x-y|^ \mu}\,dxdy.
	\end{equation*}
	\item For any $u \in W^{1,p}(\RR^N)\backslash \{0\}$, $u_t(x):=u(t^{-1}x)$ for $t>0$.
	\item The letters $C$, and $C_i$ denote various positive constants possibly different in different places.
\end{itemize}

We organize the rest of the paper as follows. In Section 3, we give some technical lemmas that will help prove our Theorem \ref{Ethm1.1}. In Section 4, we study the limit problem and prove it has a positive radial solution. In section 5, we give some background material for Theorem \ref{Ethm1.2}, and in the last section, we present its proof.

\section{Technical lemmas}
In this section, we will first collect some variational frameworks and results that form the background material. Then we prove the limiting case problem \eqref{E1.2}, using a couple of critical results.

First, we look at the Mountain-pass geometry of the functional $J_\infty$.
\begin{lem}\label{Elem3.4} The functional $J_\infty$ satisfies the following conditions:
	\begin{itemize}
		\item [(i)]There exist $\al,\rho >0$ such that $J_\infty(u) > \al$ for $0 < \|u\| \leq \rho$.
		\item [(ii)]$J_\infty(0) = 0$ and there exists $e \in W^{1,p}(\RR^N)$ with $\|e\| > \rho$ and $J_\infty(e) < 0$.
	\end{itemize}
\end{lem}
\proof(i) From \eqref{Eeq2.1} and Sobolev embedding theorem, there exist constants $C_1, C_2>0$ such that for each $u \in W^{1,p}(\RR^N)$,
$$\|F(u)\|_0 \leq C_1\|u\|^{2\cdot p_{*,\mu}}+ C_2\|u\|^{2q}.$$
Then it follows that
\begin{equation*}
	\begin{aligned}
		J_\infty(u) 
		& \geq C\left(  \|u\|^p - \|u\|^{2\cdot p_{*,\mu}}- \|u\|^{2q}\right).		\end{aligned}
\end{equation*}
It is easy to see that there exists $\al >0$ such that $J_\infty(u) > \al$, for $\|u\| \leq \rho$,  provided that $\rho$ is sufficiently small.\\
(ii) For any $u \in W^{1,p}(\RR^N)\backslash\{0\}$ and $t>0$ we define $u_t(x) := u(\frac{x}{t})$. 
\begin{equation*}
	\begin{aligned}
		J_\infty(u_t)
		= &  \frac{at^{N-p}}{p}\int\limits_{\RR^N} |\nabla u|^p dx + \frac{bt^{(N-p)\te}}{p\te}\left(\, \int\limits_{\RR^N} |\nabla u|^p dx\right)^{\te} +\frac{V_\infty t^{N}}{p}\int\limits_{\RR^N}|u|^pdx -\frac{t^{2N-\mu}}{2}\|F(u)\|_0.
	\end{aligned}
\end{equation*}
Since, $ \max\{N-p, (N-p)\te, N\} <2N-\mu $, therefore as $t \to \infty$ $J_\infty(u_t) \to -\infty$. Thus for sufficiently large $t$, we have $e = u_t \in W^{1,p}(\RR^N)$, with $\|e\| > \rho$, such that $J_\infty(e) < 0$. This completes the proof. \qed

By the classical Mountain pass theorem, we have the minimax value defined as
\begin{equation}\label{Eeq3.1}
	c_\infty := \inf\limits_{h \in \Gamma_\infty}\max\limits_{t \in [0,1]}J_\infty(h(t)),
\end{equation}
where 
\begin{equation*}
	\Gamma_\infty = \{h \in C([0,1], W^{1,p}(\RR^N)): h(0) = 0\;\text{and}\; J_\infty(h(1)) < 0\}.
\end{equation*} 
By Lemma \ref{Elem3.4}, there exists a $(PS)_{c_\infty}$ sequence $\{u_n\}$ in $W^{1,p}(\RR^N )$ at the level $c_\infty$, that is,
\begin{equation*}
	J_\infty(u_n) \to c_\infty \quad \text{and}\quad J_\infty'(u_n) \to 0 \;\text{as}\; n \to \infty.
\end{equation*}
Thus, we have established the existence of the Palais-Smale sequence. Now, to ensure that the Poho\v{z}aev identity acting on the Palais-Smale sequence approaches zero, we will use Jeanjean's technique\cite{jeanjean}.
\begin{lem}\label{Elem3.5}
	There exists a sequence $\{u_n\}$ in $W^{1,p}(\RR^N)$ such that, as $n \to \infty$,
	\begin{equation}\label{Eeq2.2}
		J_\infty(u_n) \to c_\infty, \quad  J_\infty'(u_n) \to 0 \quad \text{and}\quad P_\infty(u_n) \to 0.
	\end{equation}
	Moreover, $\{u_n\}$ is a bounded sequence in $W^{1,p}(\RR^N)$.
\end{lem}
\proof Define the map $\phi : W^{1,p}(\RR^N) \times \RR \to W^{1,p}(\RR^N)$ 
by
\begin{equation*}
	\phi(u,t) = u\left(e^{-t}x\right), \;\text{for all} \;x\in \RR^N.
\end{equation*}
For any $(u,t) \in W^{1,p}(\RR^N) \times \RR$, we introduce the composite functional  $J_\infty\circ\phi $ computed as
\begin{equation*}
	\begin{aligned}
		(J_\infty\circ\phi)(u,t) =& \frac{ae^{t(N-p)}}{p}\int\limits_{\RR^N} |\nabla u|^p dx + \frac{be^{t(N-p)\te}}{p\te}\left(\, \int\limits_{\RR^N} |\nabla u|^p dx\right)^{\te} +\frac{V_\infty e^{tN}}{p}\int\limits_{\RR^N}|u|^pdx\\ &-\frac{e^{t(2N-\mu)}}{2}\|F(u)\|_0,
	\end{aligned}
\end{equation*}
then $J_\infty\circ\phi$ is a $C^1$-functional.
Define the family of paths
\begin{equation*}
	\tilde{\Ga}_\infty = \{\tilde{\ga} \in C([0,1], W^{1,p}(\RR^N) \times \RR): \tilde{\ga}(0) = (0,0) \;\text{and}\; (J_\infty\circ\phi)(\tilde{\ga}(1)) <0\},
\end{equation*}
and the corresponding minimax level
$$ \tilde{c}_\infty= \inf\limits_{\tilde{\ga} \in \tilde{\Ga}_\infty}\max\limits_{t \in [0,1]}(J_\infty\circ\phi)(\tilde{\ga}(t)).$$
Since $\{\phi\circ\tilde{\ga} : \tilde{\ga}\in \tilde{\Ga}_\infty\} \subset \Ga_\infty $ and $\Ga_\infty \times \{0\} \subset \tilde{\Ga}_\infty $, the mountain pass levels of the functional $J_\infty$ and $J_\infty\circ\phi$ coincide, that is
\begin{equation}\label{Eeq3.2n}
	c_\infty = \tilde{c}_\infty.
\end{equation}
For each $n \in \mathbb{N}$, by \eqref{Eeq3.1}, there exists $\{g_n\} \subset \Ga_\infty$ such that $\max\limits_{t \in [0,1]}J_\infty(g_n(t)) \leq c_\infty +\frac{1}{n}$. Setting $\tilde{g}_n(t) := (g_n(t), 0) $  for $t \in [0,1]$, we see that $\tilde{g}_n\in \tilde{\Ga}_\infty$ and using \eqref{Eeq3.2n}, we get
$$\max\limits_{t \in [0,1]}(J_\infty \circ\phi)(\tilde{g}_n(t)) \leq \tilde{c}_\infty +\frac{1}{n}.$$
By Ekeland's variational principle, there exists a sequence $\{(v_n, t_n)\}\subset W^{1,p}(\RR^N)\times \RR $ such that, as $n \to \infty$,
\begin{equation}\label{Eeq3.2}
	(J_\infty\circ \phi)(v_n, t_n) \to \tilde{c}_\infty, \; (J_\infty\circ\phi)'(v_n, t_n) \to 0,
\end{equation}
$$\ds\min\limits_{t\in [0,1]}\|(v_n,t_n) - \tilde{g}_n(t)\|_{W^{1,p}(\RR^N)\times\RR} = \min\limits_{t\in [0,1]}(\|v_n-g_n\|^2 + t_n^2)^\frac{1}{2} \to 0.$$
Note that, $|t_n| \to 0$ as $n \to \infty$. Let $u_n := \phi(v_n, t_n)$, from \eqref{Eeq3.2}, we have $$J_\infty(u_n) \to c_\infty\; \text{as}\; n \to \infty.$$
It is easy to verify that for every $(w, s) \in W^{1,p}(\RR^N)\times \RR$,
\begin{equation}\label{Eeq3.3}
	\left\langle  (J_\infty\circ\phi)'(v_n,t_n), (w,s)\right\rangle  = \left\langle  J_\infty'(u_n), \phi(w,s)\right\rangle  + s P_\infty(u_n).
\end{equation}
Employing \eqref{Eeq3.3} and \eqref{Eeq3.2} with $(w,s)=(0,1)$,
we reach the conclusion that $$P_\infty(u_n)\to 0\; \text{as}\; n \to \infty.$$ 
For any $z \in W^{1,p}(\RR^N)$, set $w = \phi(z, -t_n) = z(e^{t_n}x)$ and $s= 0$ in \eqref{Eeq3.3}, we get
\begin{equation*}
	\left\langle  J_\infty'(u_n), z\right\rangle =	o(1)\|z\|,
\end{equation*}
since  $|t_n| \to 0$, as $n \to \infty$. 
That is $J_\infty'(u_n) \to 0$ in $(W^{1,p}(\RR^N))^*$. Hence, we have got a sequence $\{u_n\}$ in $W^{1,p}(\RR^N)$ that satisfies \eqref{Eeq2.2}.\\
Further, the boundedness of $(P S)_{c_\infty}$ sequence $\{u_n\}$ follows from
\begin{equation*}
	\begin{aligned}
		&J_\infty(u_n) - \frac{P_\infty(u_n)}{2N-\mu}\\
		&= \left(\frac{N-\mu+p}{2N-\mu} \right) \frac{a}{p}\int\limits_{\RR^N} |\nabla u_n |^p dx + \left(\frac{2N-\mu-\te(N-p)}{\te(2N-\mu)} \right)\frac{b}{p}\left(\, \int\limits_{\RR^N} |\nabla u_n|^p dx\right)^{\te}\\
		&\quad+\left(\frac{N-\mu}{2N-\mu} \right)\frac{V_\infty }{p}\int\limits_{\RR^N} |u_n |^p dx\\
		& \geq  C\|u_n\|^p,
	\end{aligned}
\end{equation*}
where $C >0$. This completes the proof.\qed

\begin{lem}\label{Elem3.6}
	For any $u\in W^{1,p}(\RR^N)\backslash\{0\}$, there exists a unique $\tilde{t} >0$ such that $u_{\tilde{t}} \in \mathcal{P_\infty}$, where $u_t(x)= u(t^{-1}x)$. Moreover, $J_\infty(u_{\tilde{t}})= \max\limits_{t>0}J_\infty(u_{t}). $
\end{lem}
\proof  For any $u\in W^{1,p}(\RR^N)\backslash\{0\}$ and $t>0$, we define
\begin{equation*}
	\begin{aligned}
		h_\infty(t) :=& J_\infty(u_t)\\
		=& \frac{at^{N-p}}{p}\int\limits_{\RR^N} |\nabla u|^p dx + \frac{bt^{(N-p)\te}}{p\te}\left(\, \int\limits_{\RR^N} |\nabla u|^p dx\right)^{\te} +\frac{V_\infty t^N}{p}\int\limits_{\RR^N}|u|^pdx -\frac{t^{2N-\mu}}{2}\|F(u)\|_{0}.
	\end{aligned}
\end{equation*}
Observe that for $t$ small enough, $h_\infty(t) > 0$ and $h_\infty(t) \to -\infty $ as $t \to \infty$, since $2N-\mu > \max\{N-p, (N-p)\te, N \}$. Thus, $h_\infty$ has at least one critical point, say at $\tilde{t} >0$. Further, taking the derivative of $h_\infty$, we obtain that $t\frac{d }{dt}h_\infty(t) = P_\infty(u_t)$, which implies 
$u_{\tilde{t}} \in \mathcal{P_\infty}$.\\
Now we claim that $\tilde{t}$ is unique. Let us suppose that $t_1$ and $t_2$ are two critical points of $h_\infty$, such that $0 < t_1 < t_2$.  Depending upon the range of $\te$, we consider the following cases:\\
\textbf{Case I:} $\ds1 \leq \te \leq \frac{N}{N-p}$\\
Since $P_\infty(u_{t_i}) = 0$, for $i = 1,2$ and $t_{i}>0$, we deduce that
\begin{equation*}
	\begin{aligned}
		&\frac{(N-p)a}{p}\left(\frac{1}{t_{1}^p}- \frac{1}{t_{2}^p} \right) \int\limits_{\RR^N} |\nabla u|^p dx + \frac{(N-p)b}{p}\left(\frac{1}{t_{1}^{N-(N-p)\te}} - \frac{1}{t_{2}^{N-(N-p)\te}}\right) \left(\, \int\limits_{\RR^N} |\nabla u|^p dx\right)^{\te}\\
		&= \frac{(2N-\mu)}{2}\left( t_{1}^{N-\mu}- t_{2}^{N-\mu}\right) \|F(u)\|_0.
	\end{aligned}
\end{equation*}
This is not true as $0< t_1 < t_2$.\\
\textbf{Case II:} $\ds \frac{N}{N-p} < \te < \frac{2N-\mu}{N-p}$\\
Following the same argument as in Case I, we divide by $t_{i}^{(N-p)\te}$ to get a contradiction.
This completes the proof.\qed

Set 
\begin{equation}\tag{$*$}\label{Eeq3.4}
	\mathfrak{p}_\infty = \inf\limits_{\mathcal{P_\infty}}J_\infty(u),
\end{equation} 
which is well defined, since for any $u \in \mathcal{P_\infty}$, $$\ds J_\infty(u) = 
J_\infty(u) - \frac{P_\infty(u)}{2N-\mu} \geq 0.$$

\begin{lem}\label{Elem3.7}
For $c_\infty$ as defined in \eqref{Eeq3.1} and $\mathfrak{p}_\infty$ as in \eqref{Eeq3.4}, we have $\mathfrak{p}_\infty = c_\infty.$
\end{lem}
\proof For any $\ga \in \Gamma_\infty$, we claim that $\ga([0,1]) \cap \mathcal{P_\infty} \neq \emptyset$. Indeed, there exists $\rho >0$, small enough such that for $\|u\| \leq \rho$ $$P_\infty(u)\geq 0$$ and  
\begin{equation*}
J_\infty(u) - \frac{P_\infty(u)}{2N-\mu} \geq 0, \;\;\text{for any}\; u \in W^{1,p}(\RR^N).
\end{equation*}
Hence, $\ga $ crosses the manifold $\mathcal{P_\infty}$, since 
$J_\infty(\ga(1)) < 0$. 
Therefore,
\begin{equation*}
\mathfrak{p}_\infty \leq \max\limits_{t\in [0,1]} J_\infty(\ga(t)),
\end{equation*}
which implies that $\mathfrak{p}_\infty \leq c_\infty$.\\
On the other hand, for $t_0$ large enough we see that $J_\infty(u_{t_0}) < 0$, for $u \in W^{1,p}(\RR^N)\backslash\{0\}$. Define
\begin{equation*}
\overline{\ga}(t) = \begin{cases}
u_{t_0t} & \;\text{for}\; t>0;\\
0 & \;\text{for}\; t=0,
\end{cases}
\end{equation*}
then $\overline{\ga} \in \Gamma_\infty.$ So
\begin{equation*}
c_\infty \leq \max\limits_{t \in [0,1)}J_\infty(u_{t_0t})\leq \max\limits_{t >0}J_\infty(u_{t_0t}),
\end{equation*}
which gives $c_\infty \leq \mathfrak{p}_\infty.$  \qed

\section{Radial ground states of limiting problem}
This section is devoted to the proof of Theorem \ref{Ethm1.1}. First, we will revisit the polarization technique and derive some necessary lemma to prove the radial symmetry of a given ground state solution. Then, we deploy the concentration-compactness principle to obtain the existence of a positive radial groundstate solution. 

Let $H \subset \RR^N$ be a closed half-space and that $\sigma_H$ denotes the reflection concerning its boundary $\partial H$. Define $x_H:= \sigma_H(x)$ for $x\in \RR^N$.
\begin{defi}
The polarization $u^H : \RR^N \to \RR$ of $u : \RR^N \to \RR$ is defined for $x \in \RR^N$ by
\begin{equation*}
u^H(x) = \begin{cases}
\max\{u(x), u(x_H)\} & \text{if} \; x\in H\\
\min\{u(x), u(x_H)\} & \text{if} \; x\notin H\\
\end{cases}
\end{equation*}
\end{defi}

\begin{lem}\label{Elem4.1}\cite[Lemma 5.3]{brock}
If $u \in W^{1,p}(\RR^N)$, then $u^H \in W^{1,p}(\RR^N)$. In particular, we have
\begin{equation*}
\|\nabla u\|_{L^p} = \|\nabla u^H\|_{L^p} \; \text{and}\; \|u\|_{L^p} = \|u^H\|_{L^p}.
\end{equation*}
\end{lem}
\begin{lem}
Let $u \in W^{1,p}(\RR^N)$, then 
\begin{equation}\label{Eeq4.4}
\iint\limits_{\RR^{2N}} \frac{F(u)F(u)}{|x-y|^ \mu}dxdy \leq \iint\limits_{\RR^{2N}} \frac{F(u^H)F(u^H)}{|x-y|^ \mu}dxdy.
\end{equation}
\end{lem}
\proof As $F$ is a non-decreasing function, we get that $F(u^H) = F(u)^H$. Observe that  for any $x,y \in \RR^N$ 
\begin{equation}\label{Eeq4.5}
|x-y| = |x_H-y_H| \leq |x_H- y| = |x-y_H| \text{ and } 
\end{equation}
\begin{equation}\label{Eeq4.6}
\begin{aligned}
\iint\limits_{\RR^{2N}} \frac{F(u^H(x))F(u^H(y))}{|x-y|^ \mu}dxdy =& \iint\limits_{H\times H} \frac{F(u^H(x))F(u^H(y)) + F(u^H(x_H))F(u^H(y_H))}{|x-y|^ \mu}dxdy\\ &+ \iint\limits_{H\times H} \frac{F(u^H(x_H))F(u^H(y)) + F(u^H(x))F(u^H(y_H))}{|x_H-y|^ \mu}dxdy.
\end{aligned}
\end{equation}
We divide the proof into four cases:\\
\textbf{Case I:} When $u(x_H) \leq u(x)$ and $u(y_H)\leq u(y)$\\
From the definition of polarization and \eqref{Eeq4.6}, we deduce that
\begin{equation*}
\begin{aligned}
\iint\limits_{\RR^{2N}} \frac{F(u^H(x))F(u^H(y))}{|x-y|^ \mu}dxdy =& \iint\limits_{H\times H} \frac{F(u(x))F(u(y)) + F(u(x_H))F(u(y_H))}{|x-y|^ \mu}dxdy\\ &+ \iint\limits_{H\times H} \frac{F(u(x_H))F(u(y)) + F(u(x))F(u(y_H))}{|x_H-y|^ \mu}dxdy\\
= & \iint\limits_{\RR^{2N}} \frac{F(u(x))F(u(y))}{|x-y|^ \mu}dxdy.
\end{aligned}
\end{equation*}
\textbf{Case II:} When $u(x) \leq u(x_H)$ and $u(y)\leq u(y_H)$\\
One can prove the desired result using the arguments as in Case I.\\
\textbf{Case III:} When $u(x) \leq u(x_H)$ and $u(y_H)\leq u(y)$\\
Employing the definition of polarization, \eqref{Eeq4.6} and \eqref{Eeq4.5}, we have 
\begin{equation}\label{Eeq4.7}
\begin{aligned}
&\iint\limits_{\RR^{2N}} \frac{F(u^H(x))F(u^H(y))}{|x-y|^ \mu}dxdy\\
= &  \iint\limits_{H\times H} \left[\frac{1}{|x-y|^\mu} -\frac{1}{|x_H-y|^\mu}  \right]  \left( F(u(x_H))F(u(y)) + F(u(x))F(u(y_H))\right) \;dxdy\\ &-\iint\limits_{H\times H} \left[\frac{1}{|x-y|^\mu} -\frac{1}{|x_H-y|^\mu}  \right]\left( F(u(x))F(u(y)) + F(u(x_H))F(u(y_H))\right) \;dxdy\\
& + \iint\limits_{\RR^{2N}} \frac{F(u(x))F(u(y))}{|x-y|^ \mu}dxdy.
\end{aligned}
\end{equation}
Furthermore, by direct calculation and the assumption of the case, we see that
\begin{equation*}
F(u(x_H))F(u(y)) + F(u(x))F(u(y_H)) \geq F(u(x))F(u(y)) + F(u(x_H))F(u(y_H)).
\end{equation*}
Thus, we get the desired claim from \eqref{Eeq4.7} and above inequality.\\
\textbf{Case IV:} When $u(x_H) \leq u(x)$ and $u(y)\leq u(y_H)$\\
Similar to Case III. \\Hence  the results hold. \qed
\begin{lem}\label{Elem4.3}\cite[Lemma 5.5]{Moroz1}
Let  $u \in L^{\frac{2N}{2N-\mu}}(\RR^N )$
and $H\subset \RR^N$ be a closed half-space. If $u \geq 0$, then
\begin{equation*}
\iint\limits_{\RR^{2N}} \frac{u(x)u(y)}{|x-y|^ \mu}dxdy \leq \iint\limits_{\RR^{2N}} \frac{u^H(x)u^H(y)}{|x-y|^ \mu}dxdy,
\end{equation*}
with equality if and only if either $u^H = u$ or $u^H = u \circ \sigma_H$.
\end{lem}
\begin{lem}\label{Elem4.4}\cite[Lemma 5.6]{Moroz1}
Assume that $u \in L^2(\RR^N )$ is a nonnegative. There exist $x_0 \in \RR^N$ and a non increasing function $v : (0, \infty) \to \RR$ such
that for almost every $x \in \RR^N$, $u(x) = v(|x - x_0|)$ if and only if for every closed half-space $H \subset \RR^N$, $u^H = u$ or $u^H = u\circ \sigma_H$.
\end{lem}
\begin{lem}\label{Elem4.5}
Let $u \in W^{1,p}(\RR^N)$  be a positive  ground state solution to problem \eqref{E1.2}. Then $u$ is radially symmetric to a point.
\end{lem}
\proof
Using the  strategy of  Tanaka and Jeanjean \cite{tanaka}, we can  easily  lift the  ground state solution to the path $h \in C([0,1], W^{1,p}(\RR^N))$ such  that $h(0)=0,~J_\infty(h(t))<J_\infty(u),~ J_\infty(h(1))<0$. Now applying the polarization on the path $h$, we can achieve that $u^H$ is also a critical point (for a more details, we refer to \cite{Moroz1}). Consequently, applying Lemma \ref{Elem4.1}, we deduce that
\begin{equation*}
\iint\limits_{\RR^{2N}} \frac{F(u)F(u)}{|x-y|^ \mu}dxdy \geq \iint\limits_{\RR^{2N}} \frac{F(u^H)F(u^H)}{|x-y|^ \mu}dxdy.
\end{equation*}
Employing \eqref{Eeq4.4} and Lemma \eqref{Elem4.3}, we get either $F(u^H) = F(u)$ or $F(u^H) = F(u\circ \sigma_H)$. Let us assume that $F(u^H) = F(u)$ in $\RR^N$. 
Moreover, since $F$ is non-decreasing 
\begin{equation*}
\int_{u(x_H)}^{u(x)} f(t) dt = F(u(x)) - F(u(x_H)) \begin{cases}
\geq 0 & \;\text{if} \;x \in H;\\
\leq 0 & \;\text{if}\; x \notin H;\\
\end{cases}
\end{equation*}
this implies that $u^H = u$.
Similarly in the case of $F(u^H) = F(u\circ \sigma_H)$, we can conclude that $u^H = u \circ \sigma_H$.  By Lemma \ref{Elem4.4}, we get the desired result.\qed

\textbf{ Proof of Theorem \ref{Ethm1.1}:} Let $\{u_n\}$ be a $(PS)_{c_\infty}$ sequence for $J_\infty$, by Lemma \ref{Elem3.5}, we have $\{u_n\}$ is bounded in $W^{1,p}(\RR^N )$ and $P_\infty(u_n) \to 0$ as $n \to \infty$.
Then there exists 
$A\in \RR$ such that
\begin{equation*}
\int\limits_{\RR^N}|\nabla u_n|^p \to A^p.
\end{equation*}
Define $\ds \de := \lim\limits_{n \to \infty}\sup\limits_{y \in \RR^N}\int\limits_{B_1(y)}|u_n|^p$. Now we divide the proof into the following parts:\\
\textbf{Case I:} If $\de = 0$ \\
From Lemma \ref{Elem3.2}, we see that $u_n \to 0 $ in $L^r(\RR^N)$, for all $r \in (p, p^*)$; hence $\|F(u_n)\|_0 \to 0$ as $n \to \infty.$ Also, as $P_\infty(u_n) \to 0$, we have $u_n \to 0$ in $W^{1,p}(\RR^N)$. This implies that $0 = \lim\limits_{n \to \infty} J_\infty(u_n) = c_\infty $, which is not true.\\
\textbf{Case II:} If $\de > 0$ \\
Then there exists $\{y_n\} \subset \RR^N$ such that
\begin{equation*}
\int\limits_{B_1(y_n)}|u_n|^p \geq \frac{\de}{2} >0.
\end{equation*}
Define $\overline{u}_n(x) := u_n(x+y_n)$, $\{\overline{u}_n\}$ is bounded $(PS)_{c_\infty}$ for $J_\infty$ with $P_\infty(\overline{u}_n) \to 0$ as $n \to \infty$. Then there exists $\overline{u} \in W^{1,p}(\RR^N)$ such that up to a sub-sequence 
\begin{equation}\label{Eq4.1}
\begin{aligned}
&\overline{u}_n\rightharpoonup \overline{u} \; \text{weakly in} \; W^{1,p}(\RR^N);\\
&\int\limits_{\RR^N}|\nabla \overline{u}_n|^p \to A^p,
\end{aligned}
\end{equation}
and by Fatou's Lemma, we conclude that
\begin{equation*}
\int\limits_{\RR^N}|\nabla \overline{u}|^p \leq A^p.
\end{equation*}
Let us assume that 
\begin{equation}\label{Eeq4.2}
\int\limits_{\RR^N}|\nabla \overline{u}|^p < A^p.
\end{equation}
Set
\begin{equation*}
J_{A, \infty}(u) = \frac{a+bA^{p(\te-1)}}{p}\int\limits_{\RR^N}|\nabla {u}|^p + +\frac{V_\infty }{p}\int\limits_{\RR^N}|u|^pdx -\frac{1}{2}\|F(u)\|_{0}.
\end{equation*}
Employing \eqref{Eq4.1}, we deduce that $J'_{A,\infty}(\overline{u}) = 0$, also $\overline{u}$ satisfies the Poho\v{z}aev identity associated with $J_{A, \infty}$
\begin{equation}\label{Eeq4.3}
P_{A, \infty}(u) = \left(a + bA^{p(\te-1)} \right) \frac{(N-p)}{p}\int\limits_{\RR^N} |\nabla u|^p dx + \frac{NV_\infty}{p}\int\limits_{\RR^N}|u|^pdx
-\frac{2N-\mu}{2}\|F(u)\|_0=0.
\end{equation}
Putting together \eqref{Eeq4.2} and \eqref{Eeq4.3}, we get
\begin{equation*}
0 = P_{A, \infty}(\overline{u}) > P_{\infty}(\overline{u}).
\end{equation*}
Hence by Lemma \ref{Elem3.6}, there exists $t_0 < 1$ such that $\overline{u}_{t_0} \in \mathcal{P_\infty}$. Using this information, we get
\begin{equation}\label{Eeq4.3n}
\begin{aligned}
c_\infty = \mathfrak{p}_\infty \leq& J_\infty(\overline{u}_{t_0}) - \frac{P_\infty(\overline{u}_{t_0})}{2N-\mu} \\
< & \left(\frac{N-\mu+p}{2N-\mu} \right) \frac{a}{p}\int\limits_{\RR^N} |\nabla \overline{u} |^p dx + \left(\frac{1}{\te}- \frac{N-p}{2N-\mu} \right)\frac{b}{p}\left(\, \int\limits_{\RR^N} |\nabla \overline{u}|^p dx\right)^{\te}\\
&+\left(\frac{N-\mu}{2N-\mu} \right)\frac{V_\infty}{p}\int\limits_{\RR^N} |\overline{u} |^p dx\\
\leq & \liminf\limits_{n \to \infty}\left[ J_\infty(\overline{u}_n) - \frac{P_\infty(\overline{u}_n)}{2N-\mu}\right] = c_\infty,
\end{aligned}
\end{equation}
which is not possible. So our assumption \eqref{Eeq4.2} is not true, that is
$$\ds\int\limits_{\RR^N}|\nabla \overline{u}|^p = A^p = \lim\limits_{n \to \infty}\int\limits_{\RR^N}|\nabla \overline{u}_n|^p, $$ which further implies that $J'_\infty(\overline{u}) =0$ and $\overline{u} \in \mathcal{P_\infty}.$ Following the same argument as in \eqref{Eeq4.3n} for $t_0 =1$, we get $J_\infty(\overline{u}) = c_\infty$. Thus $\overline{u}$ is a non-trivial solution of \eqref{E1.2}. Moreover, we claim that $\overline{u}$ is a ground state solution. \\
Let
\begin{equation*}
d = \inf\{J_\infty(u): u \;\text{is a solution of problem \eqref{E1.2}}\}.
\end{equation*}
Then $d \leq J_\infty(\overline{u}) = c_\infty =\mathfrak{p}_\infty \leq d$, where the last inequality holds by the definition of $\mathcal{P_\infty}$. Since
\begin{equation*}
0 = \langle J'_\infty(\overline{u}), \overline{u}^-\rangle \geq \|\overline{u}^-\|^2,
\end{equation*}
that is, $\overline{u} $ is a non-negative. Using the similar Moser's iteration as in \cite{reshmi_sarika}, we get $u \in L^\infty(\RR^N)$.  Now employing the  standard elliptic regularity and maximum principle (see \cite{pucci_max}),  we can achieve that $u$ is a positive solution. At this moment using Lemma \ref{Elem4.5} we have the required. 
%
\qed


\section{Compactness of high energy (PS)-sequences}
In this section, we will obtain all the technical results that form the background material for the proof of Theorem \ref{Ethm1.2}. We will first invoke a few elementary results associated with the Poho\v{z}aev manifold to achieve this. 
With the help of \eqref{Eeq1.5}, we define the Poho\v{z}aev manifold as follows:
\begin{equation*}
\mathcal{P} = \{u \in W^{1,p}(\RR^N)\backslash\{0\}: P(u)= 0\}.
\end{equation*}
Note that $\mathcal{P}$ contains all the critical points of $J$.

\begin{lem}\label{Elem5.2} The following holds
\begin{enumerate}
\item[(i)] There exist $\sigma_1, \sigma_2 > 0$ such that for any $u \in \mathcal{P}$, $\|u\| > \sigma_1$ and $J(u) \geq \sigma_2 >0.$
\item[(ii)] For any $u\in W^{1,p}(\RR^N)\backslash\{0\}$, there exists a unique $\tilde{t} >0$ such that $u_{\tilde{t}} \in \mathcal{P}$, where $u_t(x)= u(t^{-1}x)$. Moreover, $J(u_{\tilde{t}})= \max\limits_{t>0}J(u_{t}). $
\item[(iii)] $\{u \equiv 0\} $ is an isolated point of $P^{-1}(\{0\})$.
\item[(iv)] $\mathcal{P}$ is a closed set.
\end{enumerate}
\end{lem}
\proof \textit{(i)} By the Hardy-Littlewood-Sobolev inequality and condition $(\mathcal{V}_3)$, there exists $C>0$, such that
\begin{equation*}
\begin{aligned}
0 =P(u) &\geq  \frac{(N-p)a}{p}\int\limits_{\RR^N} |\nabla u|^p dx  +\frac{N V_\infty}{p}\int\limits_{\RR^N} |u|^pdx -C\|u\|^{2\cdot p_{*,\mu}}- C\|u\|^{2q}\\
& \geq \frac{N-p}{p}\min\{a, V_\infty\}\|u\|^p-C\|u\|^{2\cdot p_{*,\mu}}- C\|u\|^{2q}.
\end{aligned}
\end{equation*}
Since $2q > 2\cdot p_{*,\mu} > p$,  there exists $\sigma_1 > 0$ such that 
$\|u\| > \sigma_1$, for any $u \in \mathcal{P}$. Using this information, we deduce for some $C>0$ and $\sigma_2 >0$
\begin{equation*}
J(u) = 
J(u) - \frac{P(u)}{2N-\mu} \geq C\|u\|^p > \sigma_2.
\end{equation*}\\
\textit{(ii)} It follows the same analysis as in Lemma \ref{Elem3.6}. For the uniqueness, we assume that $t_1$ and $t_2$ are two critical points of function $h(t) := J(u_t) $, such that $0 < t_1 < t_2$. Depending on the range of $\te$, we consider the following cases:\\
\textbf{Case I:} $\ds1 \leq \te \leq \frac{N}{N-p}$\\
Set $\phi(t) := \ds\int\limits_{\RR^N}\left( N V(tx) + \nabla V(tx)\cdot tx\right) |u|^p\;dx$, with the help of condition $(\mathcal{V}_4)$ and some elementary calculations, we get that $\phi$ is a decreasing function. 
Since $P(u_{t_i}) = 0$, for $i = 1,2$ and $t_{i}>0$, we deduce that
\begin{equation*}
\begin{aligned}
&\frac{(N-p)a}{p}\left(\frac{1}{t_{1}^p}- \frac{1}{t_{2}^p} \right) \int\limits_{\RR^N} |\nabla u|^p dx + \frac{(N-p)b}{p}\left(\frac{1}{t_{1}^{N-(N-p)\te}} - \frac{1}{t_{2}^{N-(N-p)\te}}\right) \left(\, \int\limits_{\RR^N} |\nabla u|^p dx\right)^{\te}\\ 
& \qquad+ \frac{1}{p}\left( \phi(t_1) -\phi(t_2) \right) \\
&= \frac{(2N-\mu)}{2}\left( t_{1}^{N-\mu}- t_{2}^{N-\mu}\right) \|F(u)\|_0,
\end{aligned}
\end{equation*}
which gives us a contradiction as $0< t_1 < t_2$.\\
\textbf{Case II:} $\ds \frac{N}{N-p} < \te < \frac{2N-\mu}{N-p}$\\
Set $\ds\psi(t) := \frac{1}{t^{N(\te-1)-p\te}}\int\limits_{\RR^N}\left( N V(tx) + \nabla V(tx)\cdot tx\right) |u|^p\;dx$, using $(\mathcal{V}_3)$ and $(\mathcal{V}_4)$, we get that $\psi$ is a decreasing function. 
By similar reasoning as in Case I, we have
\begin{equation*}
\begin{aligned}
&\frac{(N-p)a}{p}\left(\frac{1}{t_{1}^{(N-p)(\te-1)}}- \frac{1}{t_{2}^{(N-p)(\te-1)}} \right) \int\limits_{\RR^N} |\nabla u|^p dx  + \frac{1}{p}\left( \psi(t_1) -\psi(t_2) \right)\\
&= \frac{(2N-\mu)}{2}\left( t_{1}^{N-\mu}- t_{2}^{N-\mu}\right) \|F(u)\|_0,
\end{aligned}
\end{equation*}
which is not true. Consequently, from both cases, we infer the uniqueness.

\textit{(iii)} From (i), we see that there exists a $\rho$ small enough such that $P(u) >0$ for $0 < \|u\| \leq \rho$.

\textit{(iv)} Since $P$ is a $C^1$ functional, $P^{-1}(\{0\}) = \mathcal{P} \cup \{0\}$ is a closed set. Now applying (iii), $P$ is a closed set. 
\qed

\begin{lem}\label{Elem5.1}
Assume  $(\mathcal{V}_1)-(\mathcal{V}_3)$ and $(\mathcal{F}_1)-(\mathcal{F}_2)$. Define  $\Ga = \{\ga \in C\left( [0,1], W^{1,p}(\RR^N)\right): \ga(0)= 0, J(\ga(1)) < 0  \}$ and $c =\min\limits_{\ga \in \Ga}\max\limits_{t \in [0,1]}J(\ga(t))$.
Then $c$ is well-defined and  $c = c_\infty$.
\end{lem}
\proof
First, we will show that $\Ga$ is non-empty. Indeed by Lemma \ref{Elem5.2}(i), there exists $\rho >0$ small enough, such that $J(u) \geq \al >0$ for $0 < \|u\| \leq \rho$. Additionally for any $u \in W^{1,p}(\RR^N)\backslash \{0\}$, $t >0$ and $\tau \in \RR^N$, we set $\ds u_{t,\tau}(x) = u\left( \frac{x-\tau}{t}\right) $, 
\begin{equation*}
\begin{aligned}
J(u_{t,\tau}) =& \frac{at^{N-p}}{p}\int\limits_{\RR^N} |\nabla u|^p dx + \frac{bt^{(N-p)\te}}{p\te}\left(\, \int\limits_{\RR^N} |\nabla u|^p dx\right)^{\te} +\frac{ t^{N}}{p}\int\limits_{\RR^N}V(xt + \tau)|u|^pdx\\ &-\frac{t^{2N-\mu}}{2}\|F(u)\|_0.
\end{aligned}
\end{equation*}
For fixed $t >0$ and $|\tau| \to \infty$, we see from $\mathcal{V}_1$
\begin{equation*}
\lim\limits_{|\tau| \to \infty}J(u_{t, \tau}) = J_\infty(u_t).
\end{equation*}
Hence for $t$ and $|\tau|$ large enough, we have $J(u_{t,\tau}) <0$, and thus $\Ga$ is non-empty and $c$ is well defined.\\
From $(\mathcal{V}_3)$, we get $J_\infty(u) \leq J(u)$, which implies that $\Ga \subset \Ga_\infty$, then
\begin{equation*}
c_\infty 
\leq \inf\limits_{h \in \Gamma}\max\limits_{t \in [0,1]}J_\infty(h(t)) \leq \inf\limits_{h \in \Gamma}\max\limits_{t \in [0,1]}J(h(t)) = c
\end{equation*}
and thus $c_\infty \leq c.$\\
To show that $c \leq c_\infty$, let $\e >0$ be arbitrary, and take 
$\ga \in \Ga_\infty$ such that $\ds\max\limits_{t\in [0,1]}J_\infty(\ga(t))\leq c_\infty +\frac{\e}{2}$. Set  $\ds\ga(t)_\tau(\cdot):= \ga(t)(\cdot -\tau)$, then $J(\ga(t)_\tau) \to J_\infty(\ga(t))$ as $|\tau| \to \infty $ for fixed $t>0$. This implies for $|\tau|$ large enough, $\ga(t)_{\tau} \in \Ga$ and 
\begin{equation*}
J(\ga(t)_\tau) \leq J_\infty(\ga(t)) +\frac{\e}{2} \;\text{for all}\; t \in[0,1].
\end{equation*}
So we get
\begin{equation*}
c \leq \max\limits_{t \in [0,1]}J(\ga(t)_\tau) \leq \max\limits_{t \in [0,1]} J_\infty(\ga(t)) +\frac{\e}{2} \leq c_\infty + \e.
\end{equation*}
This gives $c \leq c_\infty$ and hence the equality is proved. \qed

Set $$\mathfrak{p} = \inf\limits_{\mathcal{P}}J(u)>0,$$ it is well defined, by Lemma \ref{Elem5.2}(i).

\begin{lem}\label{Elem5.3}
$c= \mathfrak{p}$, where $c =\inf\limits_{\ga \in \Ga}\max\limits_{t \in [0,1]}J(\ga(t))$.
\end{lem}
\proof  The proof is similar to the proof of Lemma \ref{Elem3.7}. \qed


\begin{lem}\label{Elem5.4}
For any $u\in \mathcal{P}$, there exists a unique $t \in (0,1)$ such that $u_{t} \in \mathcal{P_\infty}$, where $u_t(x)= u(t^{-1}x)$.
\end{lem}

\proof 	Using $(\mathcal{V}_3)$, we get $P_\infty(u) < P(u)$, which implies that 
\begin{equation}\label{Eeq5.3}
P_\infty(u) < 0, \;\text{if}\; u \in  \mathcal{P}.
\end{equation}
By Lemma \ref{Elem3.6}, we see that there exists unique $\tilde{t}$ such that $P_\infty(u_{\tilde{t}}) = 0$. Hence from \eqref{Eeq5.3}, we conclude that $\tilde{t} \in (0,1)$. \qed

\begin{lem}\label{Elem5.5}
$c$ is not a critical value of the functional $J$.
\end{lem}
\proof Let us suppose that $c$ is a critical value of the functional $J$; that is, there exists $u \in W^{1,p}(\RR^N)$ such that 
$$J(u) = c, \;\; J'(u) = 0.$$
As a consequence of Lemma \ref{Elem5.3} and \ref{Elem5.4}, we see that $J(u) = \mathfrak{p}$ and $u_t \in \mathcal{P_\infty}$, for some unique $t \in (0,1)$. Further employing assumption $(\mathcal{V}_2)$, Lemma \ref{Elem3.7} and Lemma \ref{Elem5.1}, we obtain
\begin{equation*}
\begin{aligned}
\mathfrak{p} =& J(u) = J(u) - \frac{P(u)}{2N-\mu}\\
> & \left(\frac{N-\mu+p}{2N-\mu} \right) \frac{at^{N-p}}{p}\int\limits_{\RR^N} |\nabla u |^p dx + \left(\frac{2N-\mu-\te(N-p)}{\te(2N-\mu)} \right)\frac{bt^{(N-p)\te}}{p}\left(\, \int\limits_{\RR^N} |\nabla u|^p dx\right)^{\te}\\
&+\left(\frac{N-\mu}{2N-\mu} \right)\frac{V_\infty t^{N}}{p}\int\limits_{\RR^N} |u |^p dx\\
= &  J_\infty(u_t) - \frac{P_\infty(u_t)}{2N-\mu} \geq 
\mathfrak{p},
\end{aligned}
\end{equation*}
which is not possible. 
\qed

\begin{lem}\label{Elem5.6} The following holds
\begin{enumerate}
\item [(i)] $\mathcal{P}$ is a $C^1$ manifold.
\item [(ii)] Each critical point of $J\left|_\mathcal{P}\right. $ is a critical point of $J$ in $W^{1,p}(\RR^N).$
\end{enumerate}
\end{lem}
\proof \textit{(i)} We claim that $P'(u) \neq 0$, for all $u \in \mathcal{P}$. Indeed, suppose for the sake of contradiction, there exists $u \in \mathcal{P}$ such that $P'(u) = 0$, namely, $u$ is a solution to the problem
\begin{equation}\label{Eeq5.5}
\begin{aligned}
-&(N-p)\left[ a +  b\te\left(\, \int\limits_{\RR^N}|\nabla u|^p\right)^{(\te-1)}\right] \De_p u + \left( N V(x) + \nabla V(x)\cdot x\right) |u|^{p-2}u\\  
&\quad= (2N-\mu)\left( I_\al *F(u)\right)  f(u).
\end{aligned}
\end{equation}
As a consequence, $u$ satisfies the Poho\v{z}aev identity related to problem \eqref{Eeq5.5}, that is,
\begin{equation}
\begin{aligned}\label{Eeq5.6}
P_1(u) = & \frac{(N-p)^2a}{p}\int\limits_{\RR^N} |\nabla u|^p dx + \frac{(N-p)^2b\te}{p}\left(\, \int\limits_{\RR^N} |\nabla u|^p dx\right)^{\te} +\frac{N^2}{p}\int\limits_{\RR^N}V(x)|u|^pdx\\
& +\frac{2N+1}{p}\int\limits_{\RR^N}\nabla V(x)\cdot x|u|^pdx  +\frac{1}{p}\int\limits_{\RR^N}\nabla x\cdot H(x)\cdot x|u|^pdx -\frac{(2N-\mu)^2}{2}\|F(u)\|_0\\
= & 0.
\end{aligned}
\end{equation}
In view of \eqref{Eeq5.6} and $P(u) = 0$, we get
\begin{equation*}
\begin{aligned}
&\frac{(N-\mu+p)(N-p)a}{p}\int\limits_{\RR^N} |\nabla u|^p dx + \frac{[(2N-\mu)-(N-p)\te](N-p)b}{p}\left(\, \int\limits_{\RR^N} |\nabla u|^p dx\right)^{\te}\\ 
&\quad+ \frac{(N-\mu)}{p}\int\limits_{\RR^N} \left( N V(x) + \nabla V(x)\cdot x\right) |u|^pdx\\
&= \frac{N+1}{p}\int\limits_{\RR^N}\nabla V(x)\cdot x|u|^pdx  +\frac{1}{p}\int\limits_{\RR^N} x\cdot H(x)\cdot x|u|^pdx.
\end{aligned}
\end{equation*}
Using hypothesis $(\mathcal{V}_2)-(\mathcal{V}_4)$, we obtain that
\begin{equation*}
\begin{aligned}
&\frac{(N-\mu+p)(N-p)a}{p}\int\limits_{\RR^N} |\nabla u|^p dx + \frac{[(2N-\mu)-(N-p)\te](N-p)b}{p}\left(\, \int\limits_{\RR^N} |\nabla u|^p dx\right)^{\te}\\ 
&\quad+ \frac{(N-\mu)N V_\infty}{p}\int\limits_{\RR^N} |u|^pdx \; < 0,
\end{aligned}
\end{equation*}
which is absurd. Thus $\mathcal{P}$ is a $C^1$ manifold.

\textit{(ii)}  Let us suppose that $u \in \mathcal{P}$ is a critical point of the functional $J\left|_{\mathcal{P}}\right.$. Then there exists a Lagrange multiplier $\la\in \RR$ such that 
\begin{equation}\label{Eeq5.7n}
J'(u) + \la P'(u) = 0.
\end{equation}
We will show that $ \la = 0$. From \eqref{Eeq5.7n}, we see that $u$ weakly satisfies the following problem:
\begin{equation*}
\begin{aligned}
-&\left[ a(1 + \la(N-p)) + b(1+\la\te(N-p)) \left(\,\int\limits_{\RR^N}|\nabla u|^p\right)^{(\te-1)}\right]  \De_p u + (1+\la N) V(x)|u|^{p-2}u\\
&\quad+ \la\nabla V(x)\cdot x |u|^{p-2}u  = (1+ \la (2N-\mu))\left( I_\al *F(u)\right)  f(u).
\end{aligned}
\end{equation*}
Consequently, $u$ satisfies the following Poho\v{z}aev identity:
\begin{equation}\label{Eeq5.8n}
P_2(u)= P(u) + \la P_1(u) = 0,
\end{equation}
where $P_1(u)$ is defined in \eqref{Eeq5.6}. Using the fact that $P(u) = 0$ in \eqref{Eeq5.8n}, we infer that
\begin{equation*}
\begin{aligned}
&\la(N-\mu+p)\frac{a(N-p)}{p}\int\limits_{\RR^N} |\nabla u|^p dx + \la[(2N-\mu)-(N-p)\te]\frac{b(N-p)}{p}\left(\, \int\limits_{\RR^N} |\nabla u|^p dx\right)^{\te}\\ 
&\quad+ \la\frac{(N-\mu)}{p}\int\limits_{\RR^N} \left( N V(x) + \nabla V(x)\cdot x\right) |u|^pdx\\
= &\frac{\la}{p}\int\limits_{\RR^N}\nabla V(x)\cdot x|u|^pdx  +\frac{\la}{p}\int\limits_{\RR^N} \left( N\nabla V(x)\cdot x + x\cdot H(x)\cdot x\right) |u|^pdx.
\end{aligned}
\end{equation*}
From $(\mathcal{V}_2)$,$(\mathcal{V}_3)$ and $(\mathcal{V}_4)$, we see that it is only possible when $\la =0$.\qed

\begin{lem}\label{Elem5.7}
For any $u \in \mathcal{P_\infty}$,
\begin{enumerate}
\item[(i)] there exists a unique $t> 1$, such that $\ds u_t= u\left(\frac{x}{t}\right) \in \mathcal{P}$.
\item[(ii)] For any $y \in \RR^N$, there exists a unique $t_y> 1$, such that $\ds u_{y, t_y}= u\left(\frac{x-y}{t_y}\right) \in \mathcal{P}$ and
\begin{enumerate}
\item $t_y \to 1$ as $|y| \to \infty$;
\item $\sup\limits_{y\in \RR^N}t_y = T \in (1, \infty)$;
\item $t_y \to t_0$ as $y \to y_0$.
\end{enumerate}
\end{enumerate}
\end{lem}
\proof \textit{(i)}	Since $\ds P(u_t) = t \frac{d J(u_t)}{dt}$, and in view of assumption $(\mathcal{V}_3)$, we get
$$\ds P(u) =\frac{dJ(u_t)}{dt}\left|_{t=1}\right. >0, \;\text{for all}\; u \in \mathcal{P_\infty}.$$
Therefore, by the geometry of the functional $P$ and by Lemma \ref{Elem5.2} (ii), we see that there exists a unique $t >1$ such that $u_{t} \in \mathcal{P}$.

\textit{(ii)} As $u \in \mathcal{P_\infty}$, by the translation invariance of $\mathcal{P_\infty}$, we see that $u(\cdot - y) \in \mathcal{P_\infty}$, for all $y \in \RR^N$. So from (ii), there exists a unique $t_y> 1$, such that $\ds u_{y, t_y} \in \mathcal{P}$, that is,
\begin{equation}\label{Eeq5.7}
\begin{aligned}
&0= P(u_{y, t_y})\\
&= \frac{a(N-p)t_y^{N-p}}{p}\int\limits_{\RR^N} |\nabla u|^p dx + \frac{b(N-p)t_y^{(N-p)\te}}{p}\left(\, \int\limits_{\RR^N} |\nabla u|^p dx\right)^{\te} +\frac{ Nt_y^{N}}{p}\int\limits_{\RR^N}V(xt_y + y)|u|^pdx\\ 
&\quad+\frac{ t_y^{N}}{p}\int\limits_{\RR^N}\nabla V(xt_y + y)\cdot (xt_y + y)|u|^pdx-\frac{(2N-\mu)t_y^{2N-\mu}}{2}\|F(u)\|_0.
\end{aligned}
\end{equation}
\textit{(ii)(a)} Let us suppose that $t_y \to \infty$ as $|y| \to \infty$. Using $\mathcal{V}_1$ and the fact that $2N-\mu > \max\{N-p, (N-p)\te, N\}$, in \eqref{Eeq5.7}, we conclude that this supposition is not possible. So $\lim\limits_{|y|\to \infty}t_y = C \in [1, \infty). $ By \eqref{Eeq5.7} and $u \in \mathcal{P_\infty}$, we deduce that
\begin{equation*}
\begin{aligned}
&\frac{a(N-p)}{p}\left( \frac{1}{C^{N-\mu+p}}-1\right)\int\limits_{\RR^N} |\nabla u|^p dx + \frac{b(N-p)}{p}\left(\frac{1}{C^{2N-\mu-(N-p)\te}}-1 \right) \left(\, \int\limits_{\RR^N} |\nabla u|^p dx\right)^{\te}\\
&\quad+\frac{ N V_\infty}{p}\left(\frac{1}{C^{N-\mu}}-1 \right)\int\limits_{\RR^N}|u|^pdx = 0,
\end{aligned}\end{equation*}
which is only possible when $C= 1$.\\
\textit{(ii)(b)} By (a), we say that there exists $R >0$ such that for any $\e>0$, $t_y < 1+\e$, whenever $|y| > R$. We are left to check for the interval $0 \leq |y| \leq R$, let us assume that there exists $\{y_n\} \subset \RR^N$ with $0 \leq |y_n| \leq R$ such that $|t_{y_n}| \to \infty$. We get a contradiction by the same argument as in (a). Thus $\sup\limits_{y\in \RR^N}t_y = T \in (1, \infty)$.\\
\textit{(ii)(c)} From (b), we get that $\lim\limits_{y \to y_0}t_y = \tau \geq 1$. Thus from \eqref{Eeq5.7} we get that $u_{y_0, \tau} \in \mathcal{P}$ and by uniqueness $\tau = t_{y_0}$. \qed

\begin{lem}\label{Elem5.8}
Let $\{u_n\} \subset \mathcal{P}$ such that $J(u_n) \to c >0$, as $n \to \infty$, then $\{u_n\}$ is bounded in $W^{1,p}(\RR^N)$.
\end{lem}
\proof By $(\mathcal{V}_2)$ and $(\mathcal{V}_3)$, it follows that $J(u_n) = J(u_n) - \dfrac{P(u_n)}{2N-\mu} \geq C \|u_n\|^p$, for some $C> 0$. \qed

\begin{lem}\label{Elem5.9}
Let $\{u_n\}$ be a bounded $(PS)_c$ sequence for $J$ with $c > 0$,
then there exists $u\in W^{1,p}(\RR^N)$ and $A \in \RR$ such that $J'_{A}(u)=0$, where
\begin{equation*}
J_A(u) = \frac{a + bA^{p(\te -1)}}{p}\int\limits_{\RR^N}|\nabla u|^p\;dx +\frac{1}{p}\int\limits_{\RR^N}V(x)|u|^pdx -\frac{1}{2}\|F(u)\|_0,
\end{equation*}
and either 
\begin{enumerate}
\item[(i)] $u_n \to u $ in $W^{1,p}(\RR^N)$ or
\item[(ii)] there exists $l \in \NN$ and $\{y^k_n\} \subset \RR^N$ with $\{y^k_n\} \to \infty$ as $n \to \infty$ for each
$1 \leq k \leq l$, nontrivial solutions $w^1, \cdots , w^l$ of the following problem
\begin{equation*}
-(a+ bA^{p(\te-1)})\De_p(u) + V_\infty |u|^{p-2}u = \left(\,\int\limits_{\mathbb{R}^N} \frac{F(u)(y)}{|x-y|^\mu}\,dy \right) f(u),
\end{equation*}
with its related functional
\begin{equation*}
J_{A,\infty}(u) = \frac{a + bA^{p(\te -1)}}{p}\int\limits_{\RR^N}|\nabla u|^p\;dx +\frac{V_\infty}{p}\int\limits_{\RR^N}|u|^pdx -\frac{1}{2}\|F(u)\|_0,
\end{equation*}
such that
\begin{enumerate}
\item [(a)] $
c + \left[\frac{1}{p} -\frac{1}{p\te} \right]b A^{p\te} = J_A(u) + \sum\limits_{k= 1}^{l}J_{A,\infty}(w^k) $;
\item [(b)] $\|u_n-u - \sum\limits_{k =1}^{l} w^k(\cdot -y_n^k)\| \to 0$ as $n \to \infty$;
\end{enumerate}
where 
\begin{equation}\label{Eeq5.8}
\ds A^p = \int\limits_{\RR^N} |\nabla u|^p + \sum\limits_{k =1}^{l}\;\int\limits_{\RR^N} |\nabla w^k|^p,
\end{equation}
Moreover, 
\begin{equation}\label{Eeq5.9}
J_{A, \infty} (w^1) = c + \left[\frac{1}{p} -\frac{1}{p\te} \right]b A^{p(\te-1)}.
\end{equation}
\end{enumerate}
\end{lem}
\proof The proof of $(i),(ii)$ is similar to that of \cite [Proposition 3.1]{chen_liu}. It is enough to prove \eqref{Eeq5.9}. Let $\{u_n\} \subset \mathcal{P}$ be a minimizing sequence for $\mathfrak{p}$, then by Ekeland's variational principle, there exists a sequence $\{v_n\} \subset \mathcal{P}$ such that
$$J(v_n) \to \mathfrak{p}; \quad J'(v_n)_{\left|_{\mathcal{P}}\right.} \to 0; \quad \|v_n - u_n\| \to 0.$$
By Lemma \ref{Elem5.6} (ii), Lemma \ref{Elem5.3} and Lemma \ref{Elem5.8}, we conclude that $\{v_n\}$ is a bounded $(PS)_c$. Then either $(i)$ or $(ii)$ holds, but by Lemma \ref{Elem5.5}, as $c$ is not the critical value of the functional $J$, so $(i)$ cannot hold. As $\{v_n\}$ is a bounded $(PS)_c$ sequence, thus there exists a function $v \in W^{1,p}(\RR^N)$ such that $v_n \rightharpoonup v$ weakly in $W^{1,p}(\RR^N)$. We make the following cases:\\
\textbf{Case I: } When $v \neq 0$\\
We see that $J'_A(v) =0$, and $v$ satisfies the following Poho\v{z}aev identity:
\begin{equation*}
\begin{aligned}
0 = P_A(v)= &\left( a+ bA^{p(\te-1)}\right)\frac{(N-p)  }{p}\int\limits_{\RR^N} |\nabla v|^p dx +\frac{N}{p}\int\limits_{\RR^N}  V(x)|v|^pdx\\
& +\frac{1}{p}\int\limits_{\RR^N}\nabla V(x)\cdot x |v|^pdx
-\frac{2N-\mu}{2}\|F(v)\|_0 .
\end{aligned}	
\end{equation*}
Also from \eqref{Eeq5.8}, we get that 
\begin{equation}\label{Eeq5.10}
0 = P_A(v) \geq P(v).
\end{equation}
From Lemma \ref{Elem5.2}(ii) and \eqref{Eeq5.10}, there exists a unique $t \in (0,1]$ such that $v_t \in \mathcal{P}$. Further using \eqref{Eeq5.8}, we obtain that
\begin{equation}\label{Eeq5.11}
\left[1-\frac{N-p}{2N-\mu} \right] \frac{b}{p}A^{p(\te-1)} \geq \left[\frac{1}{\te}-\frac{N-p}{2N-\mu} \right] \frac{bt^{(N-p)\te}}{p} \left(\, \int\limits_{\RR^N} |\nabla v|^p\;dx\right)^{\te-1} + b\left[ \frac{1}{p} - \frac{1}{p\te}\right]A^{p(\te-1)}.
\end{equation}  
Define
\begin{equation*}
g(t) :=  \frac{t^N}{p(2N-\mu)}\int\limits_{\RR^N}\left[(N-\mu)  V(tx) - \nabla V(tx)\cdot tx\right]  |v|^p\;dx,
\end{equation*}
and taking the derivative of $g$, we get
\begin{equation*}
\begin{aligned}
g'(t) =& \frac{t^{(N-1)}(N-\mu)}{p(2N-\mu)}\int\limits_{\RR^N}\left(N V(tx) + \nabla V(tx)\cdot tx \right)|v|^p\;dx\\
&-  \frac{t^{(N-1)}}{p(2N-\mu)}\int\limits_{\RR^N}\left[ \nabla V(tx)\cdot tx 
-\left(N V(tx)\cdot tx + tx\cdot H(tx) \cdot tx \right)\right] |v|^p\;dx\geq 0.	\end{aligned}
\end{equation*}
The last inequality is attained by using the $(\mathcal{V}_2)$, $(\mathcal{V}_3)$ and $(\mathcal{V}_4)$, and it follows that for all $0 < t \leq 1$, $g(t) \leq g(1)$.
Along with this observation, we use \eqref{Eeq5.10} and \eqref{Eeq5.11}, to obtain
\begin{equation}\label{Eeq5.12}
\begin{aligned}
J_A(v) =& J_A(v) - \frac{P_A(v)}{2N-\mu}\\
\geq&  \left[1 - \frac{N-p}{2N-\mu}\right]\frac{at^{N-p}}{p}\int\limits_{\RR^N}|\nabla v|^p\;dx   + \left[\frac{1}{\te} - \frac{N-p}{2N-\mu}\right]\frac{bt^{(N-p)\te}}{p}\int\limits_{\RR^N} |\nabla v|^p \\
&+ \left[ \frac{1}{p} - \frac{1}{p\te}\right]bA^{p(\te-1)}\int\limits_{\RR^N} |\nabla v|^p + g(t)\\
= & J(v_t)- \frac{P(v_t)}{2N-\mu}+ \left[ \frac{1}{p} - \frac{1}{p\te}\right]bA^{p(\te-1)}\int\limits_{\RR^N} |\nabla v|^p\\
= & c + \left[ \frac{1}{p} - \frac{1}{p\te}\right]bA^{p(\te-1)}\int\limits_{\RR^N} |\nabla v|^p.
\end{aligned}
\end{equation}
Since $J'_{A, \infty}(w^i) = 0$, for all $i = 1, \cdots l$, arguing as in \eqref{Eeq5.12}, we claim that
\begin{equation}\label{Eeq5.14}
J_{A,\infty}(w^i) \geq c + \left[ \frac{1}{p} - \frac{1}{p\te}\right]bA^{p(\te-1)}\int\limits_{\RR^N} |\nabla w^i|^p.
\end{equation}
Using (ii)(a), \eqref{Eeq5.12}, \eqref{Eeq5.14} and \eqref{Eeq5.8}, we get that
\begin{equation}\label{Eeq5.15}
\begin{aligned}
c + \left[\frac{1}{p} -\frac{1}{p\te} \right]b A^{p\te} 
\geq& (1+l)c + \left[ \frac{1}{p} - \frac{1}{p\te}\right]bA^{p(\te-1)}\left(\, \int\limits_{\RR^N} |\nabla v|^p + \sum_{k= 1}^{l} \int\limits_{\RR^N} |\nabla w^i|^p\right)\\
\geq & 2c +  \left[ \frac{1}{p} - \frac{1}{p\te}\right]bA^{p\te}
\end{aligned}
\end{equation}
which is a contradiction. Hence $v \equiv 0 $.
Again from (ii)(a), we get that
\begin{equation}\label{Eeq5.16}
\begin{aligned}
c + \left[\frac{1}{p} -\frac{1}{p\te} \right]b A^{p\te} 
\geq & lc + \left[ \frac{1}{p} - \frac{1}{p\te}\right]bA^{p\te}.
\end{aligned}
\end{equation}
This implies $l= 1$ and thus we have proved \eqref{Eeq5.9}.\qed

\begin{lem}\label{Elem5.10}
Let $\{u_n\}\subset\mathcal{P}$ be a $(P S)_\tau$ sequence for $J$ with $\tau\in (c, 2c)$, then there exists $u\in W^{1,p}(\RR^N)\backslash\{0\}$ such that up to a sub sequence, $u_n \to u $ in $W^{1,p}(\RR^N)$.
\end{lem}
\proof By Lemma \ref{Elem5.8}, $\{u_n\}$ is bounded in $W^{1,p}(\RR^N)$. Employing Lemma \ref{Elem5.9}, we see that either
(i) or (ii) holds. We consider the following two cases\\
\textbf{Case I:} When $u \equiv 0$\\
Then from Lemma \ref{Elem5.9}, we see that (i) cannot hold since $\tau >0$. Thus (ii) holds true with $u \equiv 0$, from \eqref{Eeq5.16}, we get that
\begin{equation*}
\tau + \left[\frac{1}{p} -\frac{1}{p\te} \right]b A^{p\te} = \sum_{k= 1}^{l}J_{A,\infty}(w^k)
\geq  lc + \left[ \frac{1}{p} - \frac{1}{p\te}\right]bA^{p\te}.
\end{equation*}
Since $\tau < 2c$, we get $l =1$, when combined with \eqref{Eeq5.9} gives us a contradiction as $\tau > c$.\\
\textbf{Case II:} When $u \not\equiv 0$\\
Let us assume that (ii) holds, then by the same argument as in \eqref{Eeq5.15}, we have
\begin{equation*}
\tau + \left[\frac{1}{p} -\frac{1}{p\te} \right]b A^{p\te} = J_A(v) + \sum_{k= 1}^{l}J_{A,\infty}(w^k)
\geq  2c +  \left[ \frac{1}{p} - \frac{1}{p\te}\right]bA^{p\te},
\end{equation*}
which is not possible. Thus there exists $u\in W^{1,p}(\RR^N)\backslash\{0\}$ such that $u_n \to u $ in $W^{1,p}(\RR^N)$.\qed

\section{Proof of Theorem \ref{Ethm1.2}}
Here our target is to prove the existence of a positive solution to problem \eqref{E1.1}. More precisely, we show that there exists a critical point of the functional $J $ of high energy. For that purpose, we will use the linking theorem with the Barycenter mapping  restricted to the Poho\v{z}aev manifold $\mathcal{P}$.

\begin{defi}
Let $S$ be a closed subset of a Banach space $X$, and $Q \subset X$ be a submanifold with relative boundary $\partial Q$. We say that $S$ and $\partial Q$ link if:
\begin{enumerate}
\item[(i)] $S\cap \partial Q = \emptyset$;
\item[(ii)] for any $h \in C(X, X)$,	such that $h\left|_{\partial Q}\right. = id$, then $h(Q) \cap S \neq \emptyset$.
\end{enumerate}
\end{defi}

\begin{thm}{(\cite{linking}, Linking theorem) }
Let $S$ be a closed subset of a Banach space $X$ and $Q \subset X$ be a submanifold with relative boundary $\partial Q$. Suppose that $\phi \in C^1(X, \RR)$, satisfies the $(P S)$ condition and
\begin{enumerate}
\item[(i)] $S$ and $\partial Q$ link;
\item[(ii)] $\inf\limits_{u \in S}\phi(u) > \sup\limits_{u \in \partial Q}\phi(u);$
\item[(iii)] $\sup\limits_{u \in Q}\phi(u) < \infty$.
\end{enumerate}
Then
$$\tau = \inf\limits_{h \in H}\sup\limits_{	u\in Q}\phi(h(u)),$$
where $H = \{ h \in C(X,X) : h\left|_{\partial Q}\right. = id \}$.
Then $\tau$ is a critical value of $\phi$ and $\tau \geq \inf\limits_{u \in S}\phi(u)$.
\end{thm} 
Next, we introduce the barycenter function and its properties, which is going to be crucial for proving the existence of a solution of problem \eqref{E1.1}.
\begin{defi}\cite{barycenter}
The barycenter map $\ba \in W^{1,p}(\RR^N) \backslash\{0\} \to \RR^N$, is defined as 
$$ \ba(u) = \frac{1}{\|u\|_{L^1}}\int\limits_{\RR^N} x\tilde{u}(x)\;dx,$$
where $\ds\tilde{u}(x) = \left( \mu(u)(x) - \frac{1}{2}\max\limits_{x\in \RR^N}\mu(u)(x)\right)^+$ and $\ds\mu(u)(x) = \frac{1}{|B_1|}\int\limits_{B_1(x)} |u(y)|$. 
\end{defi}
It follows that  $\mu(u) \in L^{\infty}(\RR^N)$ is a continuous function and $\tilde{u} \in C_0(\RR^N)$. As $\tilde{u}$ has a compact support, barycenter map $\ba(u)$ is well defined and has the following properties:
\begin{enumerate}
\item[(i)] $\ba$ is continuous function in $W^{1,p}(\RR^N)\backslash\{0\}$;
\item[(ii)] If $u$ is radial, then $\ba(u) = 0$;
\item[(iii)] $\ba(u(\cdot -y)) = \ba(u) + y$, for any $y \in \RR^N;$
\item[(iv)] If $u$ is radial, then for any $y \in \RR^N$ and $t > 0$, $\ba(u(\frac{\cdot - y}{t})) = y.$
\end{enumerate}
Set 
\begin{equation*}
\mathcal{\tilde{P}} := \{u \in \mathcal{P} : \ba(u) = 0\}.
\end{equation*}
Clearly, $\mathcal{\tilde{P}}$ is non-empty.  Indeed, let $\overline{u}$ be the positive radial ground state solution of problem \eqref{E1.2}, where $\overline{u} \in \mathcal{P_\infty}$. Further, using Lemma \ref{Elem5.7} (i), we get that $\overline{u}_t \in \mathcal{P}$, for some $t >1$, and by using by Barycenter map property (ii), we deduce that $\ba(\overline{u}_t) = 0.$ Hence
$$\tilde{\mathfrak{p}} = \inf\limits_{u \in \mathcal{\tilde{P}}}J(u),$$
is well defined, it is easy ot see  that $\tilde{\mathfrak{p}} \geq c$. Moreover, we have the following result:

\begin{lem}\label{Elem6.1}
$\tilde{\mathfrak{p}} > c$.
\end{lem}
\proof We assume that $ c= p = \tilde{\mathfrak{p}}$. Let $\{u_n\} \subset \mathcal{\tilde{P}}$ be a minimizing sequence for $\tilde{\mathfrak{p}}$.  Using the same assertions as in  Lemma \ref{Elem5.9}, there exists a bounded $(P S)_c$ sequence $\{v_n\} \subset \mathcal{P}$ such that $\|v_n - u_n\| \to 0$ and from the continuity of function $\ba$, we have 
\begin{equation}\label{Eeq6.1}
\ba(v_n) = o(1), \;\text{as}\; n \to \infty.
\end{equation}
Also from \eqref{Eeq5.16}, we get  $l =1$, which subsequently gives
\begin{equation}\label{Eeq6.2}
\|v_n - w^1(\cdot -y_n^1)\| \to 0 \; \text{in}\; W^{1,p}(\RR^N)
\end{equation}
where $y_n^1 \in \RR^N$ with $|y_n^1| \to \infty$. Thus for $\ba(w^1)$ fixed, using \eqref{Eeq6.1} and \eqref{Eeq6.2}, we  obtain that
\begin{equation*}
0 = \lim\limits_{n \to \infty}|\ba(u_n)| = \lim\limits_{n \to \infty}|\ba(v_n)| = \lim\limits_{n \to \infty}|\ba(w^1) +y_n^1| = \infty,
\end{equation*}
which is impossible.  Therefore we must have $\tilde{\mathfrak{p}} > c$.\qed

Let us consider again the positive, radially symmetric, ground state solution $\overline{u}$ of problem \eqref{E1.2}. We define a continuous operator $K: \RR^N \to \mathcal{P}$ as
$$\ds K(y)(x) := \overline{u}_{y,t_y}(x)=\overline{u}\left(\frac{x-y}{t_y}\right),$$
using Lemma \ref{Elem5.7}(ii), we see that the map is well defined for all $y$ and some unique $t_y >1$. Further, using the barycenter property (iv), we get that
$$\ba(K(y)) = y.$$ \qed

\begin{lem}
Assume that there exists a constant $T > 1$, such that the assumption $(\mathcal{V}_5)$ holds true, that is,
\begin{equation*}
\sup\limits_{x \in\RR^N} V(x) \leq V_\infty + \frac{p J_\infty(\overline{u})}{T^N \int\limits_{\RR^N}|\overline{u}|^p},
\end{equation*}  
where $T = \sup\limits_{y\in \RR^N}t_y $, then $J(K(y)) < 2c$.
\end{lem}
\proof Since $J_\infty$ is translation invariant and from Lemma \ref{Elem3.6}, we see that $J_\infty(\overline{u})= \max\limits_{t>0}J_\infty(\overline{u}_{t}). $ Then from $(\mathcal{V}_5)$, we have
\begin{equation}\label{Eeq6.3}
\begin{aligned}
J(K(y)) &= J_\infty(K(y)) + \frac{1}{p}\int\limits_{\RR^{N}}\left( V(x) - V_\infty \right)|K(y)|^p\;dx\\
& \leq J_\infty(\overline{u}_{y,t_y} ) + \frac{t_y^N c}{T^N} < 2c, \end{aligned}
\end{equation}
where $c$ is as defined in Lemma \ref{Elem5.1}.\qed

\textbf{ Proof of Theorem \ref{Ethm1.2}:}
From $(\mathcal{V}_1)$ and Lemma \ref{Elem5.1}, it follows that
\begin{equation*}
\lim\limits_{|y| \to \infty} J(K(y)) = J_\infty(\overline{u}) = c,
\end{equation*}
by Lemma \ref{Elem6.1}, $c <\tilde{p}$, then there exists $R_0> 0$ such that 
\begin{equation}\label{Eeq6.4}
\max\limits_{|y| = R_0} J(K(y)) < \tilde{\mathfrak{p}}.
\end{equation}
To apply the linking theorem, we take $Q := K(\overline{B_{R_0}(0)})$ and $S := \mathcal{\tilde{P}}$, we will prove the hypothesis holds true. First, we show that $\partial Q$ and $S$ link. For any $u \in \partial Q$, we see that $\ba(u) =  \ba(K(y)) =y  $, for any $y \in \RR^N$, with $|y| = R_0 >0$, which implies that $u \not\in S$, that is, $\partial Q \cap S = \emptyset$. Next, we show that $h(Q) \cap S \neq \emptyset$, for any $h \in \mathcal{H}$, where
$$\mathcal{H }= \{h \in C(Q, \mathcal{P}): h\left|_{\partial Q} \right. = id\}.$$ We define a map $G: \overline{B_{R_0}(0)} \to \RR^N$ as
\begin{equation*}
G(y) = (\ba \circ h \circ K)(y).
\end{equation*}
The function $G$, being a composition of continuous maps, is continuous. Moreover, for any $|y| = R_0$, we have $K(y) \in \partial Q$, thus
$G(y) = y$. By the  Brouwer's fixed point theorem, there exists $y_1 \in B_{R_0} (0)$ such that $G(y_1) = 0$, which implies that $h(K(y_1)) \in S$ and thus $h(K(y_1)) \in h(Q) \cap S$. So $S$ and $\partial Q$ link. We attain the conditions of the linking theorem from \eqref{Eeq6.3} and \eqref{Eeq6.4}. Thus, there exists a $(P S)_\tau$ sequence $\{u_n\}$ for the functional $J$,
where critical level $\tau = \inf\limits_{h \in H}\sup\limits_{u\in Q}J(h(u))\geq \inf\limits_{u \in S}J(u) = \tilde{\mathfrak{p}} > c$. Whereas by the definition of $\tau$ and using \eqref{Eeq6.3}, we get that $\tau \leq \max\limits_{u\in Q}J(u) <2c$. So $\{u_n\}$ is a $(P S)_\tau$ sequence for $J$ with $\tau \in (c, 2c)$. Lemmas \ref{Elem5.8} and \ref{Elem5.10}, guarantee the existence of a nontrivial solution $u \in W^{1,p}(\RR^N)\backslash{0}$ of problem \eqref{E1.1}. 
Reasoning as in the proof of Theorem \ref{Ethm1.2}, we conclude that $u$ is a positive solution using the maximum principle. Hence the proof of the theorem is completed.  
\qed

\textbf{Acknowledgement.}
The SERB SRG grant sponsors the first author's research- SRG/ 2022/ 001946, while the second author's research is  supported by a grant from UGC (India)  with  JRF grant number: June 18-414344.


\begin{thebibliography}{9}

\bibitem{alves} C. O. Alves, F. J. S. A. Corr{\^e}a and T. F. Ma, {\it Positive solutions for a quasilinear elliptic equation of Kirchhoff type}, Computers \& Mathematics with Applications, 49 (2005), no. 1, 85-93.



\bibitem{azzollini} A. Azzollini, {\it The elliptic Kirchhoff equation in $\RR^N$ perturbed by a local nonlinearity,} Differential Integral Equations, 25 (2012), 543–554.

\bibitem{badiale} M. Badiale and E. Nabana, {\it A note on radiality of solutions of p-Laplacian equation}, Applicable analysis, 52 (1994), no. 1-4, 35-43.


\bibitem{linking} P. Bartolo, V. Benci and D. Fortunato, {\it Abstract critical point theorems and applications to some nonlinear problems with strong resonance at infinity}, Nonlinear Analysis, 7 (1983), 981–1012.

\bibitem{barycenter} T. Bartsch and T. Weth, {\it Three nodal solutions of singularly perturbed elliptic equations on domains with topology}, Annales de l'Institut Henri Poincar{\'e} C, Analyse non lin{\'e}aire, 22 (2005), 259–281.

\bibitem{reshmi_sarika} R. Biswas, S. Goyal and K. Sreenadh, {\it Multiplicity results for p-Kirchhoff modified Schr\"odinger equations with Stein-Weiss type critical nonlinearity in $R^N$}, Differential Integral Equations, 36 (2023), 247-288. 


\bibitem{brock} F. Brock and A. Solynin, {\it An approach to symmetrization via polarization}, Transactions of the American Mathematical Society, 352 (2000), no. 4, 1759-1796.



\bibitem{Chen_li} S. J. Chen and L. Li, {\it Multiple solutions for the nonhomogeneous Kirchhoff equation on $\RR^N$}, Nonlinear Analysis: Real World Applications, 14 (2013), no. 3, 1477-1486.

\bibitem{chen_liu} P. Chen, X. C. Liu, {\it Ground states for Kirchhoff equation with Hartree-type nonlinearity}, Journal of Mathematical Analysis and Applications, 473 (2019), 587-608.

\bibitem{chen_zhang_tang} S. Chen, B. Zhang, and X. Tang, {\it Existence and non-existence results for Kirchhoff-type problems with convolution nonlinearity}, Advances in Nonlinear Analysis, 9 (2020), no. 1, 148-167.


\bibitem{ding} W. Y. Ding and W. M. Ni, {\it On the existence of positive entire solutions of a semilinear elliptic equation}, Archive for Rational Mechanics and Analysis, 91 (1986), no. 4, 283-308.




\bibitem{goel_rawat_sreenadh} D. Goel, S. Rawat and K. Sreenadh, {\it Critical growth fractional Kirchhoff elliptic problems}, to appear in Advances in Differential Equations.

\bibitem{guo} Z. Guo, {\it  Ground states for Kirchhoff equations without compact condition}, Journal of Differential Equations, 259 (2015), 2884–2902.

\bibitem{he_zou} X. M. He and W. M. Zou, {\it Existence and concentration behavior of positive solutions for a Kirchhoff equation in $\RR^3$},
Journal of Differential Equations, 252 (2012), 1813-1834.

\bibitem{jeanjean} L. Jeanjean, {\it Existence of solutions with prescribed norm for semilinear elliptic equations}, Nonlinear Analysis, Theory, Methods and Applications, 28 (1997), 1633–1659.

\bibitem{jeanjean-mon} L. Jeanjean, {\it On the existence of bounded Palais–Smale sequence and application to a Landesman-Lazer type problem set on $\RR^N$}, Proceedings of the Royal Society of Edinburgh Section A, 129 (1999), 787–809. 

\bibitem{tanaka} L. Jeanjean and K. Tanaka {\it A remark on least energy solutions in $R^N$}, Proceedings of the American Mathematical Society, 131 (2003), no. 8, 2399-2408.










\bibitem{li} Y. Li, F. Li, and J. Shi, {\it Existence of a positive solution to Kirchhoff type problems without compactness conditions}, Journal of Differential Equations, 253 (2012), no. 7, 2285-2294.

\bibitem{liandni} Y. Li and W. M. Ni, {\it Radial symmetry of positive solutions of nonlinear elliptic equations in $\RR^N$}, Communications in partial differential equations, 18 (1993), no. 5-6, 1043-1054.



\bibitem{li_ye_subcritical} G. Li and H. Ye, {\it Existence of positive ground state solutions for the nonlinear Kirchhoff type equations in 
$\RR^3$}, Journal of Differential Equations, 257 (2014), no. 2, 566–600.


\bibitem{leib} E. H. Lieb and  M. Loss,   Analysis, volume 14 of graduate studies in mathematics, American Mathematical Society, Providence, Rhode Island, 2001. 

\bibitem{lu_hartree_NL} D. L\"u, {\it A note on Kirchhoff-type equations with Hartree-type nonlinearities}, Nonlinear Analysis: Theory, Methods \& Applications, 99 (2014), 35--48.

\bibitem{Moroz1} V. Moroz and J. Van Schaftingen, {\it Existence of groundstates for a class of nonlinear Choquard equations,} Transactions of American Mathematical Society, 367 (2015), 6557-6579.



\bibitem{naimen} D. Naimen, {\it The critical problem of Kirchhoff type elliptic equations in dimension four}, Journal of Differential Equations, 257 (2014), no. 4, 1168-1193.


\bibitem{pucci_max} P. Pucci and J. Serrin, {\it The maximum principle}, 73,
(2007), Springer Science \& Business Media.



\bibitem{rabino} P. H. Rabinowitz, {\it On a class of nonlinear Schr\"odinger equations}, Zeitschrift Angewandte Mathematik und Physik, 43 (1992), no. 2, 270-291.

\bibitem{rawat_sreenadh} S. Rawat and K. Sreenadh, {\it Multiple positive solutions for degenerate Kirchhoff equations with singular and Choquard nonlinearity}, Mathematical Methods in the Applied Sciences, 44 (2021), 13812- 13832.

\bibitem{sun_tang} J. J. Sun and C. L. Tang, {\it Existence and multiplicity of solutions for Kirchhoff type equations}, Nonlinear Analysis: Theory, Methods \& Applications, 74 (2011), no. 4, 1212-1222.

\bibitem{tang_cheng} X. H. Tang and B. Cheng, {\it Ground state sign-changing solutions for Kirchhoff type problems in bounded domains}, Journal of Differential Equations, 261 (2016), no. 4, 2384-2402.

\bibitem{tang_chen} X. H. Tang and S. Chen, {\it Ground state solutions of Nehari–Pohozaev type for Kirchhoff-type problems with general potentials}, Calculus of Variations and Partial Differential Equations, 56 (2017), 110.

\bibitem{wu} X. Wu, {\it Existence of nontrivial solutions and high energy solutions for Schr\"odinger-Kirchhoff-type equations in $\RR^N$},
Nonlinear Analysis: Real World Applications, 12 (2011), 1278-1287.


\bibitem{hongyu} H. Ye, {\it Positive high energy solution for Kirchhoff equation in $\RR^{3}$ with superlinear nonlinearities via Nehari-Poho\v{z}aev manifold}, Discrete \& Continuous Dynamical Systems, 35, (2015), no. 8, 3857.

\bibitem{zhang_xie_jiang} Q. Zhang, H. Xie and Y. R. Jiang, {\it Ground state solutions of Poho\v{z}aev type for Kirchhoff-type problems with general convolution nonlinearity and variable potential}, Mathematical Methods in the Applied Sciences, 2022, 1-23.


\bibitem{zhu_cao} X. Zhu and D. Cao, {\it The concentration-compactness principle in nonlinear elliptic equations}, Acta Mathematica Scientia, 9 (1989), 307–323.


\end{thebibliography}
\end{document}